\documentclass[10pt,a4paper]{article}

\usepackage[T1]{fontenc}
\usepackage[utf8]{inputenc}
\usepackage{lmodern}
\usepackage{microtype}
\usepackage{xcolor}

\usepackage[
  a4paper,
  top=1.00in,
  bottom=1.00in,
  left=1.25in,
  right=1.25in
]{geometry}
\usepackage{titlesec}

\titleformat{\section}
  {\normalfont\Large\bfseries\raggedright}
  {\thesection}{1em}{}
\titleformat{\subsection}
  {\normalfont\large\bfseries\raggedright}
  {\thesubsection}{1em}{}
\titleformat{\subsubsection}
  {\normalfont\normalsize\bfseries\raggedright}
  {\thesubsubsection}{1em}{}

\usepackage{amsmath}
\usepackage{amssymb}
\usepackage{amsthm}
\usepackage{mathtools}
\usepackage{bm}

\numberwithin{equation}{section}
\allowdisplaybreaks[1]

\theoremstyle{plain}
\newtheorem{theorem}{Theorem}[section]
\newtheorem{proposition}[theorem]{Proposition}

\theoremstyle{definition}

\theoremstyle{remark}
\newtheorem{remark}[theorem]{Remark}

\usepackage{graphicx}
\graphicspath{{figures/}}
\usepackage{subcaption}
\usepackage{booktabs}
\usepackage{multirow}
\usepackage{array}
\usepackage{siunitx}
\usepackage{algorithm}
\usepackage{algpseudocode}

\usepackage{enumitem}
\setlist{nosep}
\usepackage[numbers,sort&compress]{natbib}

\usepackage{xcolor}
\usepackage[hidelinks]{hyperref}
\usepackage[nameinlink,capitalise,noabbrev]{cleveref}
\crefname{equation}{Eq.}{Eqs.}
\Crefname{equation}{Eq.}{Eqs.}
\crefname{figure}{Fig.}{Figs.}
\Crefname{figure}{Fig.}{Figs.}

\newcommand{\Methodname}{Multiscale Kernel Framework}

\hypersetup{
  pdftitle={Operator-Adapted Multiscale Kernels for High-Contrast Elliptic Inverse Problems},
  pdfauthor={Siqing Liu; Eric Chung},
  pdfsubject={Sparse-observation recovery of multiscale states and family-restricted effective coefficients in heterogeneous elliptic inverse problems},
  pdfkeywords={Inverse elliptic problems, Optimal recovery, High-contrast media, Multiscale kernel methods, GMsFEM}
}

\title{Operator-Adapted Multiscale Kernels for High-Contrast Elliptic Inverse Problems}
\author{
Siqing Liu,
Eric Chung\thanks{Corresponding author. Email: \texttt{eric.t.chung@cuhk.edu.hk}}\\[0.5em]
\small Department of Mathematics, The Chinese University of Hong Kong,
Shatin, N.T., Hong Kong
}
\date{}

\begin{document}

\maketitle
\vspace{-2.5em}

\begin{abstract}
We consider the joint recovery of the state and diffusion coefficient in a high-contrast elliptic equation from sparse, noisy state measurements. We seek the coefficient in a parameterized family defined relative to a prescribed reference coefficient. We construct a generalized multiscale finite element basis from the reference coefficient and use it to define a finite-rank kernel for the state. A weak Galerkin residual enforces the equation without pointwise differentiation across material interfaces. For each trial coefficient, the state is obtained by solving a linear least-squares problem. Variable projection eliminates the state coordinates from the objective, and the reduced gradient is evaluated without an adjoint PDE solve.

We split the state reconstruction error into multiscale approximation, penalty relaxation, and parameter estimation terms. Under an admissible patch decomposition, the squared multiscale approximation error is bounded by a local spectral tail estimate with a contrast-independent interpolation constant. For a fixed coefficient and discretization, the penalty relaxation error is of order $O(\varepsilon^2)$ when the residual weight is proportional to $\varepsilon^{-2}$. The parameter contribution is bounded by local Lipschitz continuity of the projected state. In 50 synthetic tests based on a reference medium of contrast $10^4$, the multiscale space with the weak residual gives mean relative state and log-coefficient errors of 2.37\% and 4.43\%, respectively, the lowest among the four configurations tested.
\end{abstract}

\noindent\textbf{Keywords:} Inverse elliptic problems, Optimal recovery, High-contrast media, Multiscale kernel methods, GMsFEM

\section{Introduction}
\label{sec:introduction}

We consider the recovery of a diffusion coefficient and the corresponding state in an elliptic equation from sparse, noisy observations of the state \citep{engl1996regularization,stuart2010inverse}. We assume that the coefficient belongs to a finite-dimensional family defined relative to a reference field. In the numerical experiments, the reference field prescribes the locations of channels and barriers, and the unknowns are 20 parameters that modulate its amplitude.

One application is hydraulic tomography, where permeability is inferred from pressure measurements \citep{yeh2000hydraulic}.

Discontinuous diffusion coefficients can produce jumps in the state gradient across material interfaces. Smooth trial functions can require many degrees of freedom to approximate these gradients \citep{wendland2004scattered,babuska2000can}. Strong-form collocation poses a separate difficulty. If a trial state has a jump in its normal flux, its distributional residual contains an interface term that cannot be evaluated pointwise. A weak formulation incorporates the flux through integrals involving first derivatives. A multiscale basis constructed from the reference field provides trial functions adapted to the prescribed heterogeneity.

GMsFEM uses local spectral problems to select basis functions that represent features such as connected high-conductivity channels \citep{efendiev2013gmsfem,chung2016adaptive}. As systematically detailed in the monograph by Chung et al. \citep{chung2023multiscale}, multiscale model reduction frameworks provide an important mathematical foundation for resolving the high-contrast problems. In inverse problems, multiscale spectral spaces have been used to estimate coarse effective parameters \citep{chung2019multiscaleinversion} and to reduce the forward model \citep{cohen2023highcontrast}. We use a fixed reference basis throughout the inversion. Local residual indicators and online enrichment provide ways to adapt such bases \citep{chung2014adaptive,chung2015residual}.

Kernel and Gaussian-process formulations have been used for forward and inverse PDE problems \citep{chen2021nonlinear,cockayne2019bayesian,mora2025gp}. Here, the kernel is defined by the multiscale basis. We combine sparse observations with a weak Galerkin residual and use variable projection for the resulting separable least-squares problem \citep{golub1973differentiation,golub2003separable,olearyrust2013variable}. The weak residual uses first derivatives of the trial and test functions \citep{kharazmi2021hpvpinn,scholz2025weak}.

Our contributions are as follows.
\begin{enumerate}
    \item We construct a finite-rank state kernel from a conforming GMsFEM basis and couple it to a weak Galerkin residual for joint state and coefficient estimation.
    \item We eliminate the state coordinates by variable projection. The projected state is computed by augmented least squares, and the reduced gradient uses derivatives of the Galerkin matrix with respect to the coefficient parameters.
    \item We bound the state reconstruction error by separating the effects of multiscale approximation, penalty relaxation, and parameter estimation. The multiscale estimate depends on the local spectral tail and an admissible patch decomposition.
    \item We compare two trial spaces and two residual formulations on 50 synthetic inverse problems with a common high-contrast reference geometry. The comparison measures the effects of the state representation and the residual formulation under the stated numerical settings.
\end{enumerate}

Section~\ref{sec:problem} formulates the inverse problem. Section~\ref{sec:method} describes the multiscale space and the algorithm, Section~\ref{sec:theory} gives the error estimates, and Section~\ref{sec:experiments} presents the numerical results.

\section{Problem formulation}
\label{sec:problem}

This section establishes the continuous functional setting of the inverse coefficient problem, physically motivated by steady subsurface flow in saturated porous formations \citep{yeh2000hydraulic}. In \cref{subsec:governing-inverse-problem}, we formulate the governing boundary value problem and the sparse observation model. In \cref{subsec:ideal-recovery}, we state the ideal constrained recovery benchmark in abstract function spaces and identify the fundamental bottlenecks that necessitate multiscale kernel reduction and variational relaxation.

\subsection{Governing PDE and continuous observation model}
\label{subsec:governing-inverse-problem}

Let $D\subset\mathbb R^d$ ($d\in\{2,3\}$) be a bounded, connected Lipschitz domain with boundary $\partial D$. The forward physics is governed by the second-order, divergence-form elliptic boundary value problem:
\begin{equation}
  -\nabla\cdot\left(\kappa^{\rm phys}\nabla u^{\rm phys}\right)=f
  \quad\text{in }D,
  \qquad
  u^{\rm phys}=0
  \quad\text{on }\partial D,
  \label{eq:elliptic-pde}
\end{equation}
where $u^{\rm phys}\in H_0^1(D)$ denotes the true state field and $\kappa^{\rm phys}\in L^\infty(D)$ represents the unknown heterogeneous diffusion coefficient. The deterministic source term satisfies $f\in H^{-1}(D)$. Throughout this work, the coefficient is assumed to satisfy the strict uniform ellipticity condition:
\begin{equation}
  0<\kappa_{\min}\leq\kappa^{\rm phys}(x)\leq\kappa_{\max}<\infty
  \quad\text{a.e. in }D,
  \label{eq:elliptic-bounds}
\end{equation}
where the contrast ratio $\kappa_{\max}/\kappa_{\min}$ may be arbitrarily large, characterizing high-contrast media without scale separation.

The variational formulation of \eqref{eq:elliptic-pde} seeks $u^{\rm phys}\in H_0^1(D)$ such that
\begin{equation}
  a_{\kappa^{\rm phys}}(u^{\rm phys},v)=\ell(v)
  \quad\text{for all }v\in H_0^1(D),
  \label{eq:weak-form-continuous}
\end{equation}
where the parameter-dependent continuous bilinear form $a_\kappa: H_0^1(D)\times H_0^1(D)\to\mathbb R$ and the linear load functional $\ell: H_0^1(D)\to\mathbb R$ are defined respectively by
\begin{equation*}
  a_\kappa(u,v):=\int_D\kappa(x)\nabla u(x)\cdot\nabla v(x)\,dx,
  \qquad
  \ell(v):=\langle f,v\rangle_{H^{-1},H_0^1}.
\end{equation*}
Under the coercivity guaranteed by \eqref{eq:elliptic-bounds}, the Lax--Milgram lemma ensures the existence of a unique weak solution $u^{\rm phys}\in H_0^1(D)$.

In practical diagnostic configurations, direct field measurements are limited to sparse localized evaluations. Let $\{\ell_m^{\rm obs}\}_{m=1}^{N_{\rm obs}}\subset (H_0^1(D))^*$ denote a collection of linearly independent, bounded linear observation functionals representing localized spatial averages (mollifications). While idealized point evaluations are not bounded functionals on $H_0^1(D)$ in dimensions $d\in\{2,3\}$, they are well-defined once restricted to the continuous finite-element spaces employed in discrete implementations (see \cref{sec:method}). Gathering these functionals into a vector-valued observation operator,
\begin{equation}
  \mathcal O(u)
  :=
  \begin{bmatrix}
    \ell_1^{\rm obs}(u)&\cdots&\ell_{N_{\rm obs}}^{\rm obs}(u)
  \end{bmatrix}^{\!\top}
  \in\mathbb R^{N_{\rm obs}},
  \label{eq:obs-operator}
\end{equation}
the measured physical dataset is modeled by
\begin{equation}
  y=\mathcal O(u^{\rm phys})+\eta,
  \label{eq:observation-channel}
\end{equation}
where $\eta\in\mathbb R^{N_{\rm obs}}$ represents additive, unobservable measurement noise. Recovering the pair $(u^{\rm phys},\kappa^{\rm phys})$ exclusively from the finite vector $y$ constitutes an ill-posed inverse problem because the observations constrain only finitely many state quantities \citep{engl1996regularization,stuart2010inverse}.

\subsection{Ideal function-space recovery formulation}
\label{subsec:ideal-recovery}

To establish a principled benchmark for state and parameter estimation, we cast the joint reconstruction into the framework of optimal recovery in abstract functional spaces \citep{micchellirivlin1977survey,binev2024incomplete}. Let $\mathcal U\subseteq H_0^1(D)$ be an admissible Hilbert state space endowed with inner product $\langle\cdot,\cdot\rangle_{\mathcal U}$ and induced norm $\lVert\cdot\rVert_{\mathcal U}$, which encodes prior state regularity. Let $\mathcal V\subseteq H_0^1(D)$ denote a dense test space associated with the weak residual statements. Furthermore, let $\mathfrak K\subset L^\infty(D)$ designate the admissible set of diffusion coefficients satisfying the uniform bounds \eqref{eq:elliptic-bounds}, and let $\mathcal R:\mathfrak K\to[0,\infty]$ denote an abstract regularization functional that penalizes unphysical spatial features or incorporates prior structural knowledge.

The continuous ideal recovery problem is formulated as the following hard-constrained variational minimization:
\begin{equation}
  \begin{aligned}
    \underset{u\in\mathcal U,\,\kappa\in\mathfrak K}{\operatorname{minimize}}
    \quad&
    \frac{\lambda_u}{2}\lVert u\rVert_{\mathcal U}^2
    +\mathcal R(\kappa)
    \\
    \text{subject to}\quad&
    \mathcal O(u)=y,
    \\
    &a_\kappa(u,v)=\ell(v)
    \qquad\text{for every }v\in\mathcal V,
  \end{aligned}
  \label{eq:continuous-fixed-space-constrained-recovery}
\end{equation}
where $\lambda_u>0$ is a prescribed weighting parameter balancing the energetic norm of the state against parameter regularization.

Two issues arise in solving \eqref{eq:continuous-fixed-space-constrained-recovery}. First, noisy observations may be inconsistent with the PDE constraint, so we penalize the observation misfit. Second, smooth trial spaces may need many degrees of freedom to approximate the gradient jumps described in \cref{sec:introduction}. In \cref{sec:method}, we use a reference-based multiscale trial space and replace the exact PDE constraint by a weak residual penalty.

\section{\Methodname}
\label{sec:method}

We construct the state space from a reference coefficient and combine the data misfit with a weak Galerkin residual. We then eliminate the state coordinates by variable projection.

\subsection{Multiscale state kernel and parameter domain}
\label{subsec:multiscale-state-kernel}
\label{subsec:representation-layer}
\label{subsec:coefficient-model}

\begin{figure}[htbp]
  \centering
  \includegraphics[width=0.88\linewidth]{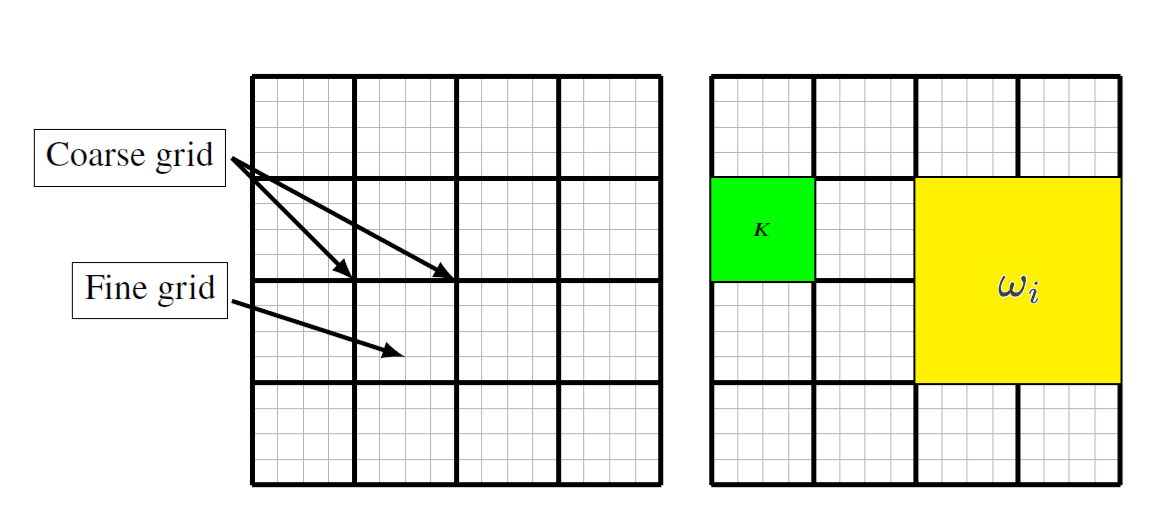}
  \caption{Coarse and fine meshes and a representative GMsFEM neighborhood.
  Thick black edges delimit \(\mathcal T_H\), thin gray lines indicate
  \(\mathcal T_h\), the green region marks one
  \(K\in\mathcal T_H\), and the yellow region is \(\omega_i\), the union of
  coarse elements sharing the central coarse vertex \(x_i\).}
  \label{fig:gmsfem-coarse-fine-neighborhood}
\end{figure}

Let $D \subset \mathbb{R}^d$ ($d \in \{2,3\}$) be a bounded polygonal or polyhedral domain. To resolve fine-scale heterogeneities while maintaining coarse-scale computational efficiency, we introduce two nested conforming spatial discretizations: a coarse mesh $\mathcal{T}_H$ with characteristic size $H$, and a fine triangulation $\mathcal{T}_h$ with mesh size $h \ll H$ that refines $\mathcal{T}_H$ and fully resolves the geometric subgrid features. Let $\mathcal{N}_H$ denote the coarse vertex index set. For each coarse vertex $x_i$, the corresponding coarse neighborhood (or patch) $\omega_i \subset D$ is defined by the union of coarse elements sharing $x_i$:
\begin{equation*}
  \omega_i := \bigcup \{K \in \mathcal{T}_H : x_i \in \overline{K}\}.
\end{equation*}
\Cref{fig:gmsfem-coarse-fine-neighborhood} illustrates the nested coarse--fine discretizations and a representative coarse patch $\omega_i$.

To capture localized high-contrast multiscale behavior without scale separation, let $\kappa_{\mathrm{ref}} \in L^\infty(D)$ denote a uniformly positive reference diffusion coefficient reflecting the background coefficient geometry. On each coarse neighborhood $\omega_i$, let $V_{\mathrm{snap}}^{(i)} \subset H^1(\omega_i)$ be a conforming fine-scale snapshot space \citep{efendiev2013gmsfem,chung2016adaptive}. Dominant local multiscale modes are extracted by solving a localized generalized eigenvalue problem: find $(\lambda_{i,k}, \psi_{i,k}) \in \mathbb{R} \times V_{\mathrm{snap}}^{(i)}$ such that
\begin{equation}
  a_i^{\mathrm{ref}}(\psi_{i,k}, v) = \lambda_{i,k} s_i^+(\psi_{i,k}, v)
  \qquad \forall v \in V_{\mathrm{snap}}^{(i)},
  \label{eq:local-spectral-problem}
\end{equation}
where the local symmetric bilinear forms are given by
\begin{align*}
  a_i^{\mathrm{ref}}(v,w) &:= \int_{\omega_i} \kappa_{\mathrm{ref}}(x) \nabla v(x) \cdot \nabla w(x) \, dx, \\
  s_i^+(v,w) &:= \int_{\omega_i} \left( 1 + \kappa_{\mathrm{ref}}(x) \sum_{j \in \mathcal{N}_H} |\nabla \chi_j^{\mathrm{ms}}(x)|^2 \right) v(x) w(x) \, dx,
\end{align*}
with $\{\chi_i^{\mathrm{ms}}\}_{i \in \mathcal{N}_H}$ being an initial multiscale partition of unity subordinate to the open cover $\{\omega_i\}$. Arranging the eigenvalues in nondecreasing order, $0 \le \lambda_{i,1} \le \lambda_{i,2} \le \cdots$, we retain the first $L_i$ dominant eigenfunctions associated with the asymptotically small eigenvalues, which capture high-permeability channels and localized inclusions \citep{chung2016adaptive}.
The auxiliary bilinear form $s_i^+$ is engineered to control local partition-of-unity gradient perturbations, ensuring that the local squared energy spectral truncation error is bounded by an inverse eigenvalue factor $\mathcal{O}(\lambda_{i,L_i+1}^{-1})$ with the stability constant $C_{\times,i}^{\mathrm{ref}}$ of the partition-of-unity nodal interpolation operator $\mathcal{M}_i$ strictly independent of the physical contrast ratio, while global approximation guarantees remain conditional on the existence of an admissible interior decomposition and coefficient comparability.

To establish global conformity within $H_0^1(D)$, these local spectral modes are assembled using fine-grid nodal interpolation $\mathcal{I}_h$ weighted by the partition of unity:
\begin{equation*}
  \varphi_{i,k} := \mathcal{I}_h(\chi_i^{\mathrm{ms}} \psi_{i,k}) \in V_{h,0},
  \qquad k = 1, \dots, L_i,
\end{equation*}
extended by zero outside $\omega_i$, where $V_{h,0} \subset H_0^1(D)$ denotes the standard piecewise linear finite element space on $\mathcal{T}_h$. Let $\mathcal{N}_H^0$ be the set containing the indices of the interior coarse vertices. Eliminating modes associated with boundary coarse nodes enforces homogeneous Dirichlet conditions, yielding the conforming multiscale state space:
\begin{equation}
  V_{\mathrm{ms},0} := \operatorname{span}\{\varphi_{i,k} : i \in \mathcal{N}_H^0, \ 1 \le k \le L_i\} \subset H_0^1(D),
  \label{eq:multiscale-space-definition}
\end{equation}
with dimension $N_c := \sum_{i \in \mathcal{N}_H^0} L_i$.

Collecting the basis functions into a feature row vector $\Phi_f(x) := [\Phi_{f,1}(x), \dots, \Phi_{f,N_c}(x)] \in \mathbb{R}^{1\times N_c}$, any state $u_c \in V_{\mathrm{ms},0}$ admits the coordinate parameterization
\begin{equation*}
  u_c(x) = \Phi_f(x) c = \sum_{j=1}^{N_c} c_j \Phi_{f,j}(x),
\end{equation*}
where $c = [c_1, \dots, c_{N_c}]^\top \in \mathbb{R}^{N_c}$. 

This construction naturally induces a finite-rank, operator-adapted reproducing kernel:
\begin{equation}
  K_{\mathrm{ms}}(x, x') := \sum_{j=1}^{N_c} \Phi_{f,j}(x) \Phi_{f,j}(x') = \Phi_f(x) \Phi_f(x')^\top.
  \label{eq:multiscale-state-kernel}
\end{equation}
For linearly independent basis functions, equip their span with the inner product $\langle\Phi_f c,\Phi_f d\rangle=c^\top d$. The inner product of $\Phi_f c$ with $K_{\mathrm{ms}}(\cdot,x')$ is $\Phi_f(x')c$. Thus $\mathcal H_{K_{\mathrm{ms}}}=V_{\mathrm{ms},0}$ and
\begin{equation}
  \|u_c\|_{K_{\mathrm{ms}}}^2 = \|c\|_2^2.
  \label{eq:multiscale-rkhs-norm}
\end{equation}
We use $\lambda_u\|c\|_2^2/2$ as the state regularization term, with $\lambda_u>0$.

For the diffusion coefficient, let $\Theta \subseteq \mathbb{R}^{N_\theta}$ be a nonempty, closed parameter domain. We consider an injective parameterization $\theta \mapsto \kappa_\theta \in L^\infty(D)$ such that each realization is strictly positive and bounded, $\kappa_\theta \in L^\infty(D)$ with $\operatorname{ess\,inf}_{x \in D} \kappa_\theta(x) > 0$. When $\Theta$ is unbounded (e.g., in unconstrained log-parameterizations), uniform ellipticity bounds $0 < \kappa_{\min}(\mathcal{K}) \le \kappa_\theta(x) \le \kappa_{\max}(\mathcal{K}) < \infty$ hold uniformly across any compact parameter subset $\mathcal{K} \subset \Theta$. Over bounded parameter domains, or when restricted to the compact parameter sublevel sets induced by coercive regularization, uniform ellipticity holds globally across the active parameter domain:
\begin{equation*}
  0 < \kappa_{\min} \le \kappa_\theta(x) \le \kappa_{\max} < \infty
  \quad \text{a.e. } x \in D, \quad \forall \theta \in \Theta.
\end{equation*}
The admissible diffusion fields form the parametric coefficient family $\mathfrak{K} := \{\kappa_\theta : \theta \in \Theta\}$, giving the joint search space $V_{\mathrm{ms},0} \times \mathfrak{K}$.
Whenever the physical coefficient $\kappa^{\rm phys}$ is well-approximated within this parameterized family, there exists a ground-truth (or best-approximation) parameter vector $\theta^\dagger \in \Theta$ inducing the in-family truth $\kappa^\dagger := \kappa_{\theta^\dagger}$ and the corresponding continuous weak solution $u^\dagger := u(\theta^\dagger) \in H_0^1(D)$ satisfying $a_{\kappa^\dagger}(u^\dagger, v) = \ell(v)$ for all $v \in H_0^1(D)$.

We assume that $\theta \mapsto \kappa_\theta$ is continuously differentiable on an open neighborhood of $\Theta$. To stabilize estimation along unresolved parameter directions, we impose a continuous, coercive parameter regularizer $R_\theta : \Theta \to [0,\infty)$, defined via the abstract functional as $R_\theta(\theta) := \mathcal{R}(\kappa_\theta)$, satisfying
\begin{equation}
  R_\theta(\theta) \longrightarrow \infty \quad \text{as } \|\theta\|_2 \longrightarrow \infty \text{ within } \Theta,
  \label{eq:parameter-regularizer-coercivity}
\end{equation}
which admits a $C^1$ extension near $\Theta$. Coercivity ensures that parameter minimizing sequences remain within compact sublevel sets, guaranteeing the existence of an optimal coefficient parameter. Due to the severe ill-posedness and finite measurement dimension, the coercive parameter regularization $R_\theta(\theta)$ is used to select a unique and stable representative from the set of parameters consistent with the observations.

\subsection{Weak-form residuals and joint relaxation}
\label{subsec:information-layer}
\label{subsec:continuous-constrained-recovery}
\label{subsec:localized-weak-information}
\label{subsec:relaxed-map-objective}

To obtain a tractable discretization of the constrained problem \eqref{eq:continuous-fixed-space-constrained-recovery}, we restrict the state $u$ to the multiscale trial space $V_{\mathrm{ms},0}$ and relax both the observation constraint and the weak PDE constraint via quadratic penalties. Because $\kappa_\theta$ is discontinuous across interfaces, the strong-form operator $-\nabla\cdot(\kappa_\theta\nabla u)$ is not defined pointwise. Enforcing the PDE via its weak Galerkin residual on $V_{\mathrm{ms},0}$ requires only $H^1(D)$ regularity and naturally preserves flux continuity across interfaces in the weak sense.

The first term incorporates the sparse observational data. Restricting the bounded linear measurement functionals $\{\ell_m^{\mathrm{obs}}\}_{m=1}^{N_{\mathrm{obs}}}$ to the multiscale space $V_{\mathrm{ms},0}$ defines the discrete observation matrix $H_f \in \mathbb{R}^{N_{\mathrm{obs}} \times N_c}$ componentwise as
\begin{equation*}
  [H_f]_{mj} := \ell_m^{\mathrm{obs}}(\Phi_{f,j}),
  \qquad m = 1, \dots, N_{\mathrm{obs}}, \quad j = 1, \dots, N_c,
\end{equation*}
which simplifies to the evaluation $[H_f]_{mj} = \Phi_{f,j}(x_m^{\mathrm{obs}})$ for point-sensor arrays.

The second term enforces variational consistency with the governing elliptic boundary value problem through a Galerkin projection onto the multiscale space $V_{\mathrm{ms},0}$. In the continuous functional setting, the parameter-dependent weak operator $A_f^{\mathrm{cont}}(\theta) \in \mathbb{R}^{N_c \times N_c}$ and the load vector $b_f^{\mathrm{cont}} \in \mathbb{R}^{N_c}$ are defined by
\begin{equation}
  [A_f^{\mathrm{cont}}(\theta)]_{ij} := \int_D \kappa_\theta(x) \nabla \Phi_{f,j}(x) \cdot \nabla \Phi_{f,i}(x) \, dx,
  \qquad
  [b_f^{\mathrm{cont}}]_i := \langle f, \Phi_{f,i} \rangle_{H^{-1}, H_0^1}.
  \label{eq:continuous-localized-weak-forms}
\end{equation}
Algebraically, let $P_{\mathrm{ms}} \in \mathbb{R}^{N_h \times N_c}$ denote the prolongation matrix whose columns store the fine-grid nodal representations of $\{\Phi_{f,j}\}_{j=1}^{N_c}$. These operators admit the compact Galerkin representations $A_f^{\mathrm{cont}}(\theta) = P_{\mathrm{ms}}^\top A_h(\kappa_\theta) P_{\mathrm{ms}}$ and $b_f^{\mathrm{cont}} = P_{\mathrm{ms}}^\top b_h$, where $A_h(\kappa_\theta)$ and $b_h$ denote the standard fine-scale stiffness matrix and load vector on $\mathcal{T}_h$, respectively.

At the discrete implementation level, evaluating the continuous integrals in \eqref{eq:continuous-localized-weak-forms} employs standard numerical quadrature over the fine triangulation $\mathcal{T}_h$, yielding the assembled weak stiffness matrix $A_f(\theta) \in \mathbb{R}^{N_c \times N_c}$ and assembled load vector $b_f \in \mathbb{R}^{N_c}$. Throughout the theoretical analysis in Section~\ref{sec:theory}, we adopt the standard assumption of exact Galerkin assembly ($A_f(\theta) = A_f^{\mathrm{cont}}(\theta)$ and $b_f = b_f^{\mathrm{cont}}$) to isolate multiscale reduction and penalty errors from quadrature defects; in practical computations, fine-mesh quadrature rules are chosen sufficiently accurate to maintain coercivity and render assembly errors negligible. Given $A_f(\theta)$ and $b_f$, the discrete weak-form PDE residual vector is given by
\begin{equation}
  r_W(c, \theta) := A_f(\theta) c - b_f.
  \label{eq:implemented-weak-information-model}
\end{equation}

To reconcile differing physical units and balance relative confidence between the data-misfit and PDE residual terms, we introduce symmetric positive definite (SPD) weight matrices $W_y \in \mathbb{R}^{N_{\mathrm{obs}} \times N_{\mathrm{obs}}}$ and $W_W \in \mathbb{R}^{N_c \times N_c}$ ($W_y \succ 0$, $W_W \succ 0$), which induce the weighted Euclidean norms $\|r\|_{W_y}^2 := r^\top W_y r$ and $\|r\|_{W_W}^2 := r^\top W_W r$. Relaxing the rigid, generally incompatible equality constraints of \eqref{eq:continuous-fixed-space-constrained-recovery} via quadratic penalty terms yields the unconstrained joint regularized objective:
\begin{equation}
  J(c, \theta) := \frac{1}{2} \|H_f c - y\|_{W_y}^2 + \frac{1}{2} \|A_f(\theta)c - b_f\|_{W_W}^2 + \frac{\lambda_u}{2} \|c\|_2^2 + R_\theta(\theta).
  \label{eq:joint-regularized-objective}
\end{equation}
The variational functional \eqref{eq:joint-regularized-objective} achieves an optimal trade-off among observational fidelity, weak PDE residual consistency, multiscale RKHS regularity, and parameter plausibility.

\begin{remark}[Bayesian Interpretation and MAP Equivalence]
\label{rem:bayesian-interpretation}
Although \eqref{eq:joint-regularized-objective} is formulated as a deterministic variational objective, it possesses an exact statistical correspondence within the Bayesian inversion framework \citep{stuart2010inverse}. Suppose the observation vector and discrete weak residual are perturbed by conditionally independent, zero-mean Gaussian noise:
\begin{equation}
  y \mid c \sim \mathcal{N}\left(H_f c, W_y^{-1}\right),
  \qquad
  b_f \mid (c, \theta) \sim \mathcal{N}\left(A_f(\theta) c, W_W^{-1}\right),
  \label{eq:observation-likelihood}
\end{equation}
where the weak residual model acts as a Gaussian pseudo-likelihood for the deterministic load balance. Assigning mutually independent prior distributions to the state coordinates and parameter vector,
\begin{equation}
  c \sim \mathcal{N}\left(0, \lambda_u^{-1} I_{N_c}\right),
  \qquad
  \pi_\theta(\theta) \propto \exp\bigl(-R_\theta(\theta)\bigr),
  \label{eq:coordinate-priors}
\end{equation}
where we assume that the parameter prior is proper, satisfying $0 < \int_\Theta \exp\bigl(-R_\theta(\theta)\bigr) \, d\theta < \infty$ (which holds unconditionally for standard quadratic regularization $R_\theta(\theta) = \frac{\alpha}{2}\|\theta\|_2^2$),
Bayes' theorem implies that the joint negative log-posterior density satisfies
\begin{equation*}
  -\log \pi(c, \theta \mid y, b_f) = J(c, \theta) + \text{const}.
\end{equation*}
Consequently, every global minimizer $(c_\star, \theta_\star)$ of \eqref{eq:joint-regularized-objective} corresponds to a maximum a posteriori (MAP) point estimate of the joint posterior distribution.
\end{remark}

\subsection{Reduced formulation via variable projection}
\label{subsec:variable-projection}

The joint recovery problem governed by \eqref{eq:joint-regularized-objective} exhibits a classical \emph{separable nonlinear least-squares} structure \citep{golub1973differentiation,olearyrust2013variable}: for any fixed parameter vector $\theta \in \Theta$, the functional $c \mapsto J(c, \theta)$ is strictly quadratic with respect to the multiscale state coordinates $c \in \mathbb{R}^{N_c}$. Factoring the symmetric positive definite weight matrices via their Cholesky decompositions, $W_y = L_y^\top L_y$ and $W_W = L_W^\top L_W$, we recast the state-dependent terms in \eqref{eq:joint-regularized-objective} into an augmented linear least-squares problem:
\begin{equation*}
  D(\theta) :=
  \begin{bmatrix}
    L_y H_f \\
    L_W A_f(\theta) \\
    \sqrt{\lambda_u} I_{N_c}
  \end{bmatrix}
  \in \mathbb{R}^{(N_{\mathrm{obs}} + 2N_c) \times N_c},
  \qquad
  d :=
  \begin{bmatrix}
    L_y y \\
    L_W b_f \\
    0
  \end{bmatrix}
  \in \mathbb{R}^{N_{\mathrm{obs}} + 2N_c}.
\end{equation*}
In this augmented representation, the joint objective takes the compact form:
\begin{equation*}
  J(c, \theta) = \frac{1}{2} \|D(\theta)c - d\|_2^2 + R_\theta(\theta).
\end{equation*}
Because the state regularization parameter satisfies $\lambda_u > 0$, the normal Gram matrix obeys the uniform spectral lower bound:
\begin{equation*}
  D(\theta)^\top D(\theta)
  = H_f^\top W_y H_f + A_f(\theta)^\top W_W A_f(\theta) + \lambda_u I_{N_c}
  \succeq \lambda_u I_{N_c} \succ 0,
  \qquad \forall \theta \in \Theta.
\end{equation*}
Consequently, $D(\theta)$ maintains full column rank for every $\theta \in \Theta$, which guarantees that the conditional state recovery problem is uniformly $\lambda_u$-strongly convex, independently of the rank or conditioning of the discrete observation operator $H_f$ and the weak stiffness matrix $A_f(\theta)$.

For each parameter $\theta \in \Theta$, there exists a unique conditionally optimal state coordinate vector:
\begin{equation}
  c_\star(\theta) := \operatorname*{argmin}_{c \in \mathbb{R}^{N_c}} J(c, \theta)
  = \operatorname*{argmin}_{c \in \mathbb{R}^{N_c}} \frac{1}{2} \|D(\theta)c - d\|_2^2,
  \label{eq:unique-projected-state}
\end{equation}
which formally satisfies the symmetric positive definite normal equations:
\begin{equation}
  \left( H_f^\top W_y H_f + A_f(\theta)^\top W_W A_f(\theta) + \lambda_u I_{N_c} \right) c_\star(\theta)
  = H_f^\top W_y y + A_f(\theta)^\top W_W b_f.
  \label{eq:projected-state-normal-equations}
\end{equation}
To preserve backward numerical stability and prevent squaring the spectral condition number of $D(\theta)$, we avoid forming the normal matrix directly and instead compute $c_\star(\theta)$ via the thin QR factorization of the augmented design matrix:
\begin{equation}
  D(\theta) = Q_D(\theta) R_D(\theta),
  \qquad
  R_D(\theta) c_\star(\theta) = Q_D(\theta)^\top d,
  \label{eq:projected-state-qr}
\end{equation}
where $Q_D(\theta) \in \mathbb{R}^{(N_{\mathrm{obs}} + 2N_c) \times N_c}$ has orthonormal columns and $R_D(\theta) \in \mathbb{R}^{N_c \times N_c}$ is nonsingular and upper triangular.

By conditionally eliminating the optimal state coordinates, we eliminate $c$ from the optimization problem and define the reduced objective functional over the parameter domain:
\begin{equation*}
  \varphi(\theta) := J(c_\star(\theta), \theta),
  \qquad \theta \in \Theta.
\end{equation*}
Under the standing assumption that $\theta \mapsto A_f(\theta)$ and $R_\theta(\theta)$ are $C^1$ on an open neighborhood of $\Theta$, the state mapping $\theta \mapsto c_\star(\theta)$ and the reduced functional $\varphi(\theta)$ inherit $C^1$ regularity. The coercivity of $R_\theta$ and the closedness of $\Theta$ ensure that $\varphi$ achieves at least one global minimizer on $\Theta$.

By the multivariate chain rule, the gradient of the reduced functional satisfies
\begin{equation*}
  \nabla \varphi(\theta) = \nabla_\theta J(c_\star(\theta), \theta) + \left( D_\theta c_\star(\theta) \right)^\top \nabla_c J(c_\star(\theta), \theta).
\end{equation*}
Since $\nabla_c J(c_\star(\theta),\theta)=0$, the chain rule gives
\begin{equation*}
  \nabla \varphi(\theta) = \nabla_\theta J(c_\star(\theta), \theta).
\end{equation*}
Differentiating $J$ directly with respect to $\theta_r$ yields the explicit reduced gradient:
\begin{equation}
  \frac{\partial \varphi(\theta)}{\partial \theta_r}
  = r_W(c_\star(\theta), \theta)^\top W_W \frac{\partial A_f(\theta)}{\partial \theta_r} c_\star(\theta)
  + \frac{\partial R_\theta(\theta)}{\partial \theta_r},
  \qquad r = 1, \dots, N_\theta.
  \label{eq:reduced-objective-gradient}
\end{equation}
The reduced gradient uses the projected state, the weak residual, and the parameter derivatives of the Galerkin matrix and regularizer. Its evaluation requires no adjoint solve or derivative of $c_\star$.

\subsection{Offline--online computational architecture}
\label{subsec:offline-online-architecture}
\label{subsec:offline-online-algorithm}

In the offline stage, we construct the local snapshot spaces, solve the reference spectral problems, and assemble the conforming multiscale basis and prolongation matrix. We also precompute the load vector, the observation matrix for a fixed sensor design, and the sparsity pattern for element integration.

\begin{figure}[t!]
  \centering
  \includegraphics[width=0.97\linewidth]{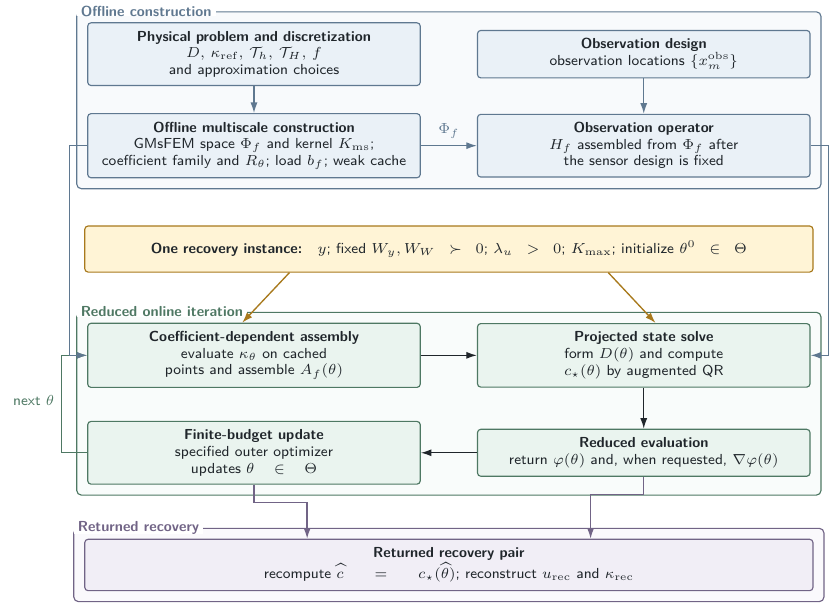}
  \caption{Offline--online workflow of the multiscale kernel framework. The reference GMsFEM state space and observation operator are constructed offline, while online evaluations invoke variable projection to solve for the conditionally linear state and update the coefficient parameter.}
  \label{fig:method-architecture}
\end{figure}

For each trial coefficient, the online stage assembles the reduced Galerkin matrix, solves the augmented least-squares problem by thin QR, and evaluates the reduced objective and gradient. The least-squares solve uses only the multiscale coordinates. The coefficient-dependent matrix entries and their derivatives can still require fine-grid quadrature. See \cref{alg:projected-state-evaluation,fig:method-architecture}.

\begin{algorithm}[htbp]
  \caption{Variable-Projection Evaluation Algorithm}
  \label{alg:projected-state-evaluation}
  \begin{algorithmic}[1]
    \Require Trial parameter $\theta \in \Theta$; precomputed structures $H_f$, $y$, $b_f$, assembly cache, weight matrices $W_y \succ 0$, $W_W \succ 0$, regularizer $R_\theta$, and state penalty parameter $\lambda_u > 0$.
    \Ensure Conditionally optimal state $c_\star(\theta)$, reduced objective value $\varphi(\theta)$, and (optionally) reduced gradient $\nabla \varphi(\theta)$.
    \State Assemble the multiscale weak stiffness matrix $A_f(\theta)$ and evaluate parameter sensitivities $\frac{\partial A_f(\theta)}{\partial \theta_r}$ using cached geometric integrals.
    \State Form the augmented design matrix $D(\theta)$ and target vector $d$.
    \State Compute the thin QR factorization $D(\theta) = Q_D(\theta) R_D(\theta)$ and solve the nonsingular upper-triangular system $R_D(\theta) c_\star(\theta) = Q_D(\theta)^\top d$ via backward substitution.
    \State Evaluate the discrete weak residual $r_W(c_\star(\theta), \theta) = A_f(\theta) c_\star(\theta) - b_f$.
    \State Evaluate the reduced objective value:
    \begin{equation*}
      \varphi(\theta) = \frac{1}{2} \|H_f c_\star(\theta) - y\|_{W_y}^2 + \frac{1}{2} \|r_W(c_\star(\theta), \theta)\|_{W_W}^2 + \frac{\lambda_u}{2} \|c_\star(\theta)\|_2^2 + R_\theta(\theta).
    \end{equation*}
    \If{the reduced gradient is requested}
      \State Compute $\nabla \varphi(\theta) \in \mathbb{R}^{N_\theta}$ componentwise via \cref{eq:reduced-objective-gradient}.
    \EndIf
    \State \Return $c_\star(\theta)$, $\varphi(\theta)$, and $\nabla \varphi(\theta)$.
  \end{algorithmic}
\end{algorithm}

Upon convergence of the outer optimizer to a terminal parameter estimate $\widehat{\theta} \in \Theta$, the conditionally optimal multiscale coordinate vector is recovered as $\widehat{c} := c_\star(\widehat{\theta})$. The physical state and heterogeneous diffusion fields are then reconstructed via
\begin{equation}
  u_{\mathrm{rec}}(x) := \Phi_f(x) \widehat{c},
  \qquad
  \kappa_{\mathrm{rec}}(x) := \kappa_{\widehat{\theta}}(x).
  \label{eq:returned-recovery-pair}
\end{equation}
To monitor observational fidelity and weak variational consistency upon termination, we evaluate the normalized relative observation mismatch and weak-form PDE residual:
\begin{equation*}
  \mathcal{R}_{\mathrm{obs}} := \frac{\|H_f \widehat{c} - y\|_2}{\|y\|_2},
  \qquad
  \mathcal{R}_{\mathrm{weak}} := \frac{\|A_f(\widehat{\theta}) \widehat{c} - b_f\|_2}{\|b_f\|_2}.
\end{equation*}
By freezing the fine-scale multiscale geometry offline and conditionally projecting out the linear state coordinates online, this architecture achieves scalable and stable parameter recovery for high-contrast media. This decoupled structural hierarchy provides the exact foundation for the comprehensive error decomposition and convergence analysis established in \cref{sec:theory}.

\section{Mathematical Analysis}
\label{sec:theory}

This section establishes the analytical foundation of the proposed multiscale
kernel recovery framework.  We first characterize the finite-dimensional
operator structure and the projected-state reduction introduced in
\cref{sec:method}.  We then establish an essential three-level state-error
decomposition that separates the multiscale kernel representation error, the
finite-penalty variational relaxation error, and the parametric optimization
error.  These results provide rigorous deterministic guarantees while
substantiating the Bayesian maximum a posteriori (MAP) perspective.

\subsection{Finite-Rank State Space and Projected-State Structure}
\label{subsec:finite-rank-projected-state}

\paragraph{Standing assumptions.}
We collect the foundational assumptions underlying the theoretical analysis:
\begin{enumerate}[label=(A\arabic*),leftmargin=2.5em]
  \item \label{ass:features}
    The multiscale basis functions $\{\Phi_{f,j}\}_{j=1}^{N_c}\subset V_{\rm ms,0}\subset H_0^1(D)$
    are linearly independent, and the observation functionals defining $H_f$
    are bounded linear operators on $H_0^1(D)$.
  \item \label{ass:coefficients}
    The reference coefficient $\kappa_{\rm ref}\in L^\infty(D)$ and the candidate
    family $\{\kappa_\theta:\theta\in\Theta\}$ satisfy the uniform ellipticity
    condition: there exist constants $0<\kappa_{\min}\le\kappa_{\max}<\infty$ such that
    $\kappa_{\min}\le\kappa_\theta(x)\le\kappa_{\max}$ almost everywhere in $D$ for all
    $\theta\in\Theta$.
  \item \label{ass:parameter-regularity}
    The admissible parameter domain $\Theta\subseteq\mathbb R^{N_\theta}$ is nonempty and
    closed. The operator map $\theta\mapsto A_f(\theta)$ and the regularizer $R_\theta$
    admit $C^1$ extensions on an open neighborhood of $\Theta$. Moreover,
    $R_\theta:\Theta\to[0,\infty)$ is continuous and coercive:
    $R_\theta(\theta)\to\infty$ as $\|\theta\|_2\to\infty$ within $\Theta$.
  \item \label{ass:metrics}
    The weight matrices satisfy $W_y\succ0$ and $W_W\succ0$, and the state
    regularization parameter is strictly positive, $\lambda_u>0$.
\end{enumerate}

\begin{proposition}[Finite-rank state space and deterministic regularization]
\label[proposition]{prop:finite-rank-joint-map}
Under the standing assumptions \ref{ass:features}--\ref{ass:metrics}, the
operator-adapted kernel defined in \cref{eq:multiscale-state-kernel} has native
reproducing kernel Hilbert space
\begin{equation}
  \mathcal H_{K_{\rm ms}}=V_{\rm ms,0},
  \qquad
  \langle u_c,u_d\rangle_{K_{\rm ms}}=c^\top d,
  \qquad
  \lVert u_c\rVert_{K_{\rm ms}}^2=\lVert c\rVert_2^2,
  \label{eq:finite-rank-rkhs-identity}
\end{equation}
for all $u_c=\Phi_f c$ and $u_d=\Phi_f d\in V_{\rm ms,0}$.  Consequently, the isotropic
coordinate penalty $\frac{\lambda_u}{2}\lVert c\rVert_2^2$ coincides identically with the
deterministic RKHS/Tikhonov functional $\frac{\lambda_u}{2}\lVert u_c\rVert_{K_{\rm ms}}^2$.

Moreover, under the conditionally independent Gaussian models established in \cref{rem:bayesian-interpretation}, any global minimizer of $J(c,\theta)$ corresponds to a joint MAP point estimate.
\end{proposition}

The proof is given in \cref{app:finite-rank-joint-map-proof}.

\begin{proposition}[Projected state and reduced-objective structure]
\label[proposition]{prop:projected-state-reduced-objective}
Let the standing assumptions \ref{ass:features}--\ref{ass:metrics} hold.  For each
fixed $\theta\in\Theta$, the state subproblem $c\mapsto J(c,\theta)$ is uniformly
$\lambda_u$-strongly convex and possesses the unique minimizer $c_\star(\theta)$
characterized by \cref{eq:unique-projected-state}.  The augmented design matrix
$D(\theta)$ satisfies the uniform spectral lower bound
\begin{equation}
  \sigma_{\min}(D(\theta))\ge\sqrt{\lambda_u}>0,
  \qquad\forall\theta\in\Theta,
  \label{eq:design-matrix-spectral-bound}
\end{equation}
and its thin QR factorization \cref{eq:projected-state-qr} stably computes $c_\star(\theta)$.

Moreover, the state mapping $\theta\mapsto c_\star(\theta)$ and the reduced objective
functional $\varphi(\theta):=J(c_\star(\theta),\theta)$ are continuously differentiable ($C^1$).
The reduced functional $\varphi$ is coercive on $\Theta$, attains at least one global
minimizer on $\Theta$, and its gradient is given explicitly by the adjoint-free formula \cref{eq:reduced-objective-gradient}.

If, additionally, $\theta\mapsto A_f(\theta)$ and $R_\theta$ are real analytic on a neighborhood
of $\Theta$, then $c_\star$ and $\varphi$ are real analytic.
\end{proposition}

The proof is given in \cref{app:projected-state-reduced-objective-proof}.

\cref{prop:projected-state-reduced-objective} guarantees that the conditional state problem is
well-posed independently of the rank or conditioning of the individual measurement
channel $H_f$ or weak operator $A_f(\theta)$.  The reduced gradient in
\cref{eq:reduced-objective-gradient} eliminates state sensitivities, rendering
online parameter optimization computationally scalable.

\subsection{State error decomposition}
\label{subsec:state-error-decomposition}

To analyze the accuracy of the reconstructed state $u_{\rm rec}$, we first split the total error $\lVert u^\dagger - u_{\rm rec}\rVert_X$ into three intermediate errors using the triangle inequality, and then bound each term individually in \cref{subsec:multiscale-approximation,subsec:relaxation-sensitivity}.

Throughout this subsection, we assume the target diffusion coefficient lies within the admissible parameter family: there exists $\theta^\dagger\in\Theta$ such that $\kappa^\dagger=\kappa_{\theta^\dagger}$. Let $u^\dagger\in H_0^1(D)$ denote the exact continuous physical state satisfying
\begin{equation}
  a_{\kappa^\dagger}(u^\dagger,v)=\langle f,v\rangle_{H^{-1},H_0^1}
  \qquad\forall v\in H_0^1(D).
  \label{eq:true-continuous-pde}
\end{equation}
We define the conforming multiscale Galerkin benchmark state $u_G(\theta^\dagger)\in V_{\rm ms,0}$ as the exact solution to the unpenalized weak Galerkin balance at the truth parameter:
\begin{equation}
  A_f(\theta^\dagger)c_G(\theta^\dagger)=b_f,
  \qquad
  u_G(\theta^\dagger):=\Phi_f c_G(\theta^\dagger).
  \label{eq:multiscale-galerkin-benchmark}
\end{equation}
By the ellipticity assumption \ref{ass:coefficients}, $A_f(\theta^\dagger)$ is symmetric positive definite, ensuring $u_G(\theta^\dagger)$ is uniquely defined.

For any parameter vector $\theta\in\Theta$, let $u_\star(\theta):=\Phi_f c_\star(\theta)$ be the conditionally projected state. Let $\theta_\star\in\arg\min_{\theta\in\Theta}\varphi(\theta)$ denote an exact global minimizer of the reduced objective, and let $\widehat\theta\in\Theta$ denote the computational parameter iterate returned by the finite-budget optimization algorithm.
The final reconstructed state is given by $u_{\rm rec}=\Phi_f c_\star(\widehat\theta)=u_\star(\widehat\theta)$.

\begin{proposition}
\label[proposition]{prop:three-level-state-error-decomposition}
Let $X$ be any normed space continuously embedded in $H_0^1(D)$ on which the states $u^\dagger$, $u_G(\theta^\dagger)$, and $u_\star(\theta)$ are well defined. Under the standing assumptions \ref{ass:features}--\ref{ass:metrics}, the total state reconstruction error satisfies:
\begin{equation}
  \lVert u^\dagger - u_{\rm rec}\rVert_X
  \le
  \underbrace{
    \lVert u^\dagger - u_G(\theta^\dagger)\rVert_X
  }_{E_{\rm ms}}
  +
  \underbrace{
    \lVert u_G(\theta^\dagger) - u_\star(\theta^\dagger)\rVert_X
  }_{E_{\rm relax}}
  +
  \underbrace{
    \lVert u_\star(\theta^\dagger) - u_\star(\widehat\theta)\rVert_X
  }_{E_{\rm opt}}.
  \label{eq:three-level-state-error-decomposition}
\end{equation}
Furthermore, the parametric term $E_{\rm opt}$ admits the natural separation into an inverse-selection defect and an outer-optimization defect:
\begin{equation}
  E_{\rm opt}
  \le
  \underbrace{
    \lVert u_\star(\theta^\dagger) - u_\star(\theta_\star)\rVert_X
  }_{E_{\rm inv}}
  +
  \underbrace{
    \lVert u_\star(\theta_\star) - u_\star(\widehat\theta)\rVert_X
  }_{E_{\rm opt}^{\rm num}}.
  \label{eq:parametric-sub-decomposition}
\end{equation}
\end{proposition}

\begin{proof}
Inserting and subtracting $u_G(\theta^\dagger)$ and $u_\star(\theta^\dagger)$ within the norm $\lVert u^\dagger - u_\star(\widehat\theta)\rVert_X$ and applying the triangle inequality yields \cref{eq:three-level-state-error-decomposition}. Applying the triangle inequality with respect to the intermediate state $u_\star(\theta_\star)$ yields \cref{eq:parametric-sub-decomposition}.
\end{proof}

The decomposition \cref{eq:three-level-state-error-decomposition} cleanly separates the distinct components of the numerical recovery:
\begin{itemize}[leftmargin=2em]
  \item $E_{\rm ms}$ reflects the approximation capacity of the operator-adapted multiscale space $V_{\rm ms,0}$, governed by the local spectral truncation error (analyzed in \cref{subsec:multiscale-approximation});
  \item $E_{\rm relax}$ captures the penalty relaxation error of the weak residual at the fixed true parameter $\theta^\dagger$ (analyzed in \cref{subsec:relaxation-sensitivity});
  \item $E_{\rm opt}$ measures the parameter-induced state error governed by the regularity of the parameter-to-state map $\theta\mapsto u_\star(\theta)$ (analyzed in \cref{subsec:relaxation-sensitivity}).
\end{itemize}

\begin{remark}[Out-of-family model mismatch]
\label[remark]{rem:model-mismatch}
If the true physical coefficient $\kappa^{\rm phys}$ does not belong to the declared family $\mathfrak K=\{\kappa_\theta:\theta\in\Theta\}$, the total error requires an initial model-mismatch projection $\inf_{\theta\in\Theta}\lVert u^{\rm phys}-u^\dagger(\theta)\rVert_X$, which precedes the three-level decomposition above. Our analysis focuses primarily on the in-family recovery regime.
\end{remark}

\subsection{Multiscale Kernel Representation Analysis}
\label{subsec:multiscale-approximation}

We now analyze the first component of the error decomposition in \cref{prop:three-level-state-error-decomposition}, specialized to the target energy norm $\lVert\cdot\rVert_{a_{\theta^\dagger}}$:
\begin{equation}
  E_{\rm ms} = \lVert u^\dagger - u_G(\theta^\dagger)\rVert_{a_{\theta^\dagger}},
  \label{eq:multiscale-representation-error-def}
\end{equation}
where $u_G(\theta^\dagger)=\Phi_f c_G(\theta^\dagger)\in V_{\rm ms,0}$ is the conforming multiscale Galerkin benchmark state defined in \cref{eq:multiscale-galerkin-benchmark}.

To isolate the spatial discretization defect from the multiscale model reduction, let $V_{h,0}\subset H_0^1(D)$ denote the conforming fine-grid piecewise linear finite element space, and let $u_h(\theta^\dagger)\in V_{h,0}$ be the fine-scale Galerkin approximation satisfying
\begin{equation}
  a_{\theta^\dagger}(u_h(\theta^\dagger), v_h) = \ell(v_h)
  \qquad \forall v_h\in V_{h,0}.
  \label{eq:fine-scale-galerkin-state}
\end{equation}
Because the reference multiscale space satisfies the conforming nesting $V_{\rm ms,0}\subset V_{h,0}\subset H_0^1(D)$, the multiscale benchmark state $u_G(\theta^\dagger)$ is equivalently the Galerkin projection of $u_h(\theta^\dagger)$ onto $V_{\rm ms,0}$ under the bilinear form $a_{\theta^\dagger}(\cdot,\cdot)$.

\begin{proposition}[Best approximation and Pythagorean Galerkin split]
\label[proposition]{prop:multiscale-best-approximation} Under the standing assumptions \ref{ass:features}--\ref{ass:metrics}, the conforming nesting $V_{\rm ms,0}\subset V_{h,0}\subset H_0^1(D)$, and assuming exact Galerkin assembly ($A_f(\theta^\dagger) = A_f^{\mathrm{cont}}(\theta^\dagger)$ and $b_f = b_f^{\mathrm{cont}}$), the multiscale representation error satisfies
\begin{equation}
  E_{\rm ms}
  =
  \inf_{v_{\rm ms}\in V_{\rm ms,0}}
  \lVert u^\dagger - v_{\rm ms}\rVert_{a_{\theta^\dagger}}.
  \label{eq:exact-multiscale-best-approximation}
\end{equation}
Moreover, it admits the exact orthogonal decomposition
\begin{equation}
  E_{\rm ms}^2
  =
  \lVert u^\dagger - u_h(\theta^\dagger)\rVert_{a_{\theta^\dagger}}^2
  +
  \lVert u_h(\theta^\dagger) - u_G(\theta^\dagger)\rVert_{a_{\theta^\dagger}}^2.
  \label{eq:continuous-fine-multiscale-split}
\end{equation}
\end{proposition}

\begin{proof}
By Galerkin orthogonality, $a_{\theta^\dagger}(u^\dagger - u_G(\theta^\dagger), v_{\rm ms}) = 0$ for all $v_{\rm ms}\in V_{\rm ms,0}$.  The Pythagorean theorem gives \cref{eq:exact-multiscale-best-approximation}.  Similarly, because $V_{\rm ms,0}\subset V_{h,0}$, the fine-grid discretization error $u^\dagger - u_h(\theta^\dagger)$ is $a_{\theta^\dagger}$-orthogonal to $V_{h,0}$, which contains $u_h(\theta^\dagger) - u_G(\theta^\dagger)$.  Decomposing $u^\dagger - u_G(\theta^\dagger) = (u^\dagger - u_h(\theta^\dagger)) + (u_h(\theta^\dagger) - u_G(\theta^\dagger))$ and expanding the energy norm yields \cref{eq:continuous-fine-multiscale-split}.
\end{proof}

The first term in \cref{eq:continuous-fine-multiscale-split} represents the fine-scale Galerkin discretization defect $\lVert u^\dagger - u_h(\theta^\dagger)\rVert_{a_{\theta^\dagger}}$, which achieves the classical rate $\mathcal O(h)$ under piecewise $H^2$-regularity on interface-fitted meshes and is made negligible by choosing $\mathcal T_h$ sufficiently fine to resolve all geometric subgrid inclusions.

The second term represents the multiscale reduction error incurred by approximating fine-grid features via the finite-rank kernel space $V_{\rm ms,0}$.  To bound this term, we decompose fine-scale functions into local representations and control the discarded spectral tails.

\paragraph{Local spectral decomposition and admissible representations.}
For each interior coarse vertex $i\in\mathcal I_0:=\mathcal N_H^0$, let $W_i:=V_{\rm snap}^{(i)}$ denote the local snapshot space on $\omega_i$.  Recall the localized generalized eigenvalue problem \cref{eq:local-spectral-problem}:
\begin{equation}
  a_i^{\rm ref}(\psi_{i,k}, v) = \lambda_{i,k} s_i^+(\psi_{i,k}, v)
  \qquad \forall v\in W_i,
  \label{eq:local-spectral-problem-recap}
\end{equation}
where the eigenvalues are ordered non-decreasingly, $0\le\lambda_{i,1}\le\lambda_{i,2}\le\cdots$.  Let $\Pi_i^{(L_i)}$ denote the $s_i^+$-orthogonal projection operator onto ${\rm span}\{\psi_{i,k}\}_{k=1}^{L_i}$, and let $\lambda_{i,L_i+1}>0$ denote the first omitted eigenvalue.

Local spectral functions are mapped into $V_{h,0}$ via the partition-weighted nodal interpolation operator:
\begin{equation}
  \mathcal M_i z := \mathcal I_h(\chi_i^{\rm ms} z) \in V_{h,0},
  \qquad z\in W_i.
  \label{eq:nodal-assembly-operator}
\end{equation}
Because boundary coarse modes are removed to guarantee zero trace on $\partial D$, any fine-grid function $w\in V_{h,0}$ must be represented exclusively using interior coarse patches.  We define the space of admissible interior decompositions of $w$ by:
\begin{equation}
  \mathcal D_h(w)
  :=
  \left\{
    \boldsymbol r = (r_i)_{i\in\mathcal I_0} : r_i\in W_i,
    \quad
    w = \sum_{i\in\mathcal I_0} \mathcal M_i r_i
  \right\}.
  \label{eq:admissible-local-decompositions}
\end{equation}
Throughout our multiscale analysis, we adopt the working hypothesis that target fine-grid states admit an admissible interior decomposition $\boldsymbol r \in \mathcal D_h(w)$ (i.e., $\mathcal D_h(w) \neq \emptyset$), which holds when boundary partition functions do not fully isolate fine nodes near coarse boundary corners (see the technical discussion in \cref{app:multiscale-approximation-proof}).

For any choice $\boldsymbol r\in\mathcal D_h(w)$, the retained spectral modes yield a conforming multiscale approximant $v_{\rm ms} := \sum_{i\in\mathcal I_0} \mathcal M_i \Pi_i^{(L_i)} r_i \in V_{\rm ms,0}$, with truncation error:
\begin{equation}
  w - v_{\rm ms} = \sum_{i\in\mathcal I_0} \mathcal M_i e_i,
  \qquad
  e_i := r_i - \Pi_i^{(L_i)} r_i.
  \label{eq:spectral-residual-decomposition}
\end{equation}

\paragraph{Spectral tail decay and contrast-independent stability.}
Define the local spectral operator $\mathcal A_i:W_i\to W_i$ by $s_i^+(\mathcal A_i z, v) = a_i^{\rm ref}(z,v)$ for all $v\in W_i$.  The energy and mass of the discarded local component $e_i$ are controlled by:
\begin{equation}
  a_i^{\rm ref}(e_i,e_i) + s_i^+(e_i,e_i)
  \le
  \frac{1}{\lambda_{i,L_i+1}} \mathfrak R_i(r_i),
  \qquad
  \mathfrak R_i(r_i) := a_i^{\rm ref}(r_i,r_i) + s_i^+(\mathcal A_i r_i, \mathcal A_i r_i).
  \label{eq:local-spectral-tail-bound}
\end{equation}
To assemble local estimates globally, we define the nodal-product stability constant on each patch:
\begin{equation}
  C_{\times,i}^{\rm ref}
  :=
  \sup_{0\neq z\in W_i}
  \frac{\lVert\mathcal M_i z\rVert_{a_{\rm ref}}^2}
       {a_i^{\rm ref}(z,z) + s_i^+(z,z)}.
  \label{eq:nodal-product-constant}
\end{equation}
Taking the infimum over all admissible decompositions defines the multiscale error functional:
\begin{equation}
  \mathcal E_h(w)
  :=
  \inf_{\boldsymbol r\in\mathcal D_h(w)}
  \sum_{i\in\mathcal I_0}
  \frac{C_{\times,i}^{\rm ref}}{\lambda_{i,L_i+1}}
  \mathfrak R_i(r_i).
  \label{eq:infimal-multiscale-functional}
\end{equation}

To transfer bounds from the frozen reference coefficient $\kappa_{\rm ref}$ to the online target parameter $\theta$, we define the essential comparability ratio:
\begin{equation}
  \rho_+(\theta) := \operatorname*{ess\,sup}_{x\in D} \frac{\kappa_\theta(x)}{\kappa_{\rm ref}(x)} < \infty.
  \label{eq:reference-comparability-ratio}
\end{equation}

\begin{theorem}
\label[theorem]{thm:implementation-matched-infimal-multiscale}
Assume a shape-regular triangulation $\mathcal T_h$ consisting of elementwise non-obtuse elements (in $d=2,3$), where coarse patches overlap with a finite maximum count $N_{\rm ov}$ (with $N_{\rm ov}\le4$ in standard two-dimensional triangulations). Assume that $\kappa_{\rm ref}$ is cellwise constant on $\mathcal T_h$, and that the initial multiscale partition functions $\{\chi_i^{\rm ms}\}$ are piecewise linear with uniformly bounded gradients. Furthermore, assume that the fine-grid state $u_h(\theta)$ admits an admissible decomposition $\boldsymbol r \in \mathcal D_h(u_h(\theta))$ within the local snapshot spaces, and that exact Galerkin assembly holds. Then there exists a geometric constant $C_{\rm sh}>0$, strictly independent of the mesh sizes $h, H$ and the physical contrast ratio $\kappa_{\max}/\kappa_{\min}$, such that
\begin{equation}
  C_{\times,i}^{\rm ref} \le C_{\rm sh}
  \qquad \forall i\in\mathcal I_0.
  \label{eq:uniform-nodal-product-bound}
\end{equation}
Furthermore, for any parameter $\theta\in\Theta$, the multiscale reduction squared energy error satisfies:
\begin{equation}
  \lVert u_h(\theta) - u_G(\theta)\rVert_{a_\theta}^2
  \le
  \rho_+(\theta) N_{\rm ov} \mathcal E_h(u_h(\theta)).
  \label{eq:implementation-matched-infimal-estimate}
\end{equation}
\end{theorem}

The proof is given in \cref{app:multiscale-approximation-proof}.

\begin{remark}[Contrast-independent robustness of the multiscale kernel]
\label[remark]{rem:contrast-robustness}
\Cref{thm:implementation-matched-infimal-multiscale} establishes that the nodal-product stability constant $C_{\rm sh}$ in \cref{eq:uniform-nodal-product-bound} remains bounded uniformly independent of the contrast ratio $\kappa_{\max}/\kappa_{\min}\to\infty$. Under a controlled decomposition numerator, the multiscale squared energy error decays with the spectral factor $\mathcal O(\lambda_{\min,L+1}^{-1})$ (corresponding to an energy-norm decay of $\mathcal O(\lambda_{\min,L+1}^{-1/2})$), while potential residual contrast dependence is isolated within the comparability ratio $\rho_+(\theta)$, the omitted eigenvalues, and the decomposition norms $\mathfrak R_i(r_i)$.
\end{remark}

\subsection{Variational Relaxation and State Sensitivity}
\label{subsec:relaxation-sensitivity}

This subsection controls the remaining two components of the three-part state error split in \cref{prop:three-level-state-error-decomposition}: the penalty relaxation error $E_{\rm relax}$ and the parameter estimation error $E_{\rm opt}$.

\paragraph{Penalty relaxation error.}
Recall the multiscale Galerkin benchmark state $u_G(\theta^\dagger)=\Phi_f c_G(\theta^\dagger)$ from \cref{eq:multiscale-galerkin-benchmark}, which satisfies the unpenalized balance $A_f(\theta^\dagger)c_G(\theta^\dagger)=b_f$.  The relaxed state $u_\star(\theta^\dagger)=\Phi_f c_\star(\theta^\dagger)$ instead incorporates the weak physical residual as a finite quadratic penalty within $J(c,\theta^\dagger)$.  To analyze how $u_\star(\theta^\dagger)$ approximates $u_G(\theta^\dagger)$ as the physical penalty confidence increases, we fix a symmetric positive definite (SPD) base weight matrix $\overline W_W\succ0$ and consider the parameterized weight family:
\begin{equation}
  W_W(\varepsilon) := \varepsilon^{-2}\overline W_W,
  \qquad \varepsilon>0.
  \label{eq:penalty-parameter-scaling}
\end{equation}
Let $c_\varepsilon(\theta)$ denote the corresponding unique projected state coordinate vector satisfying \cref{eq:projected-state-normal-equations} with weight matrix $W_W(\varepsilon)$.  We also define the $\theta$-independent normal matrix $B$ and data vector $h$:
\begin{equation}
  B := H_f^\top W_y H_f + \lambda_u I_{N_c}\succ0,
  \qquad
  h := H_f^\top W_y y.
  \label{eq:observation-regularization-block}
\end{equation}

\begin{proposition}[Finite-penalty variational relaxation]
\label[proposition]{prop:finite-penalty-relaxation}
Let the standing assumptions \ref{ass:features}--\ref{ass:metrics} hold.  For any parameter $\theta\in\Theta$, the operator $A_f(\theta)$ is nonsingular, and the relaxation displacement satisfies the exact identity:
\begin{equation}
  c_\varepsilon(\theta) - c_G(\theta)
  =
  \varepsilon^2
  \left(A_f(\theta)^\top \overline W_W A_f(\theta) + \varepsilon^2 B\right)^{-1}
  \left(h - B c_G(\theta)\right).
  \label{eq:exact-finite-penalty-displacement}
\end{equation}
Consequently, for any normed state space $X$ continuously embedded in $H_0^1(D)$,
\begin{equation}
  \lVert\Phi_f(c_\varepsilon(\theta) - c_G(\theta))\rVert_X
  \le
  \lVert\Phi_f\rVert_{\ell^2\to X}
  \frac{\varepsilon^2}
       {\sigma_{\min}\big(\overline W_W^{1/2} A_f(\theta)\big)^2}
  \lVert h - B c_G(\theta)\rVert_2.
  \label{eq:conditioning-explicit-finite-penalty-bound}
\end{equation}
In particular, representing the implemented weak metric as $W_W = \varepsilon_0^{-2}\overline W_W$, the relaxation error at the true parameter satisfies:
\begin{equation}
  E_{\rm relax}
  =
  \lVert u_G(\theta^\dagger) - u_\star(\theta^\dagger)\rVert_X
  =
  \mathcal O(\varepsilon_0^2).
  \label{eq:relaxation-error-asymptotic-rate}
\end{equation}
\end{proposition}

The proof is given in \cref{app:finite-penalty-relaxation-proof}.

\Cref{prop:finite-penalty-relaxation} demonstrates that for any fixed coefficient parameter, relaxing the equality PDE constraint into a quadratic penalty introduces a state perturbation that decays quadratically as $\varepsilon\to 0$.  The constant depends explicitly on the spectral conditioning of the weak operator and the compatibility defect $h - B c_G(\theta^\dagger)$ between the measurement channel and the unpenalized Galerkin state.

\paragraph{Projected-state sensitivity and optimization defect.}
We now bound the parametric error $E_{\rm opt} = \lVert u_\star(\theta^\dagger) - u_\star(\widehat\theta)\rVert_X$.
Let $\mathcal N\subset\Theta$ be a prescribed compact convex subset containing the true parameter $\theta^\dagger$, an exact global reduced minimizer $\theta_\star$, and the returned algorithmic iterate $\widehat\theta$ (whose existence is assumed, for instance, when $\Theta$ is convex and $\mathcal N$ is taken as the closed convex hull of a compact sublevel set of $\varphi$).

Define the uniform projected-state sensitivity constant on $\mathcal N$ by:
\begin{equation}
  L_c(\mathcal N) := \sup_{\theta\in\mathcal N} \lVert D_\theta c_\star(\theta)\rVert_{2\to2}.
  \label{eq:compact-state-lipschitz-constant}
\end{equation}

\begin{proposition}[Projected-state Lipschitz continuity and parametric reduction]
\label[proposition]{prop:projected-state-lipschitz-reductions}
Under the standing assumptions \ref{ass:features}--\ref{ass:metrics}, the sensitivity constant satisfies $L_c(\mathcal N) < \infty$.  Moreover, for all $\theta, \mu\in\mathcal N$,
\begin{equation}
  \lVert u_\star(\theta) - u_\star(\mu)\rVert_X
  \le
  \lVert\Phi_f\rVert_{\ell^2\to X} L_c(\mathcal N) \lVert\theta - \mu\rVert_2.
  \label{eq:compact-projected-state-lipschitz}
\end{equation}
In particular, the total parametric error $E_{\rm opt}$ and its two sub-components satisfy:
\begin{equation}
  \begin{aligned}
    E_{\rm opt}
    &\le
    \lVert\Phi_f\rVert_{\ell^2\to X} L_c(\mathcal N) \lVert\theta^\dagger - \widehat\theta\rVert_2, \\
    E_{\rm inv}
    &\le
    \lVert\Phi_f\rVert_{\ell^2\to X} L_c(\mathcal N) \lVert\theta^\dagger - \theta_\star\rVert_2, \\
    E_{\rm opt}^{\rm num}
    &\le
    \lVert\Phi_f\rVert_{\ell^2\to X} L_c(\mathcal N) \lVert\theta_\star - \widehat\theta\rVert_2.
  \end{aligned}
  \label{eq:inverse-optimization-state-reductions}
\end{equation}
\end{proposition}

The explicit derivative formula and proof are given in \cref{app:projected-state-lipschitz-proof}.

\begin{remark}[Scope of parameter convergence, identifiability, and conditional stability]
\label[remark]{rem:identifiability-scope}
\Cref{prop:projected-state-lipschitz-reductions} establishes that the parameter-to-state map $\theta\mapsto u_\star(\theta)$ is locally Lipschitz continuous, transferring any parameter error directly to the reconstructed state.  We emphasize that unconditionally bounding the underlying parameter distances $\lVert\theta^\dagger - \theta_\star\rVert_2$ and $\lVert\theta_\star - \widehat\theta\rVert_2$ requires establishing coefficient identifiability under sparse measurements and global convergence guarantees for non-convex optimization algorithms, both of which remain challenging open problems in PDE-constrained inverse problems \citep{engl1996regularization,stuart2010inverse}.

Nevertheless, conditional deterministic stability can be established under standard local growth conditions: if the reduced objective functional $\varphi(\theta)$ satisfies a local quadratic growth condition with modulus $\mu > 0$ on a neighborhood of a strict local minimizer $\theta_\star$, i.e.,
\begin{equation}
  \varphi(\theta) - \varphi(\theta_\star) \ge \frac{\mu}{2} \lVert\theta - \theta_\star\rVert_2^2,
  \label{eq:local-quadratic-growth}
\end{equation}
then any numerical estimate $\widehat\theta$ that achieves the objective suboptimality tolerance $\varphi(\widehat\theta) - \varphi(\theta_\star) \le \varepsilon_{\rm tol}$ and lies within the growth neighborhood $U(\theta_\star)$ satisfies the explicit parametric error bound:
\begin{equation}
  \lVert\theta_\star - \widehat\theta\rVert_2 \le \sqrt{\frac{2\varepsilon_{\rm tol}}{\mu}}.
  \label{eq:param-error-bound-growth}
\end{equation}
Combining \cref{eq:param-error-bound-growth} with the state sensitivity estimate in \cref{eq:inverse-optimization-state-reductions}, the algorithmic state error is explicitly bounded by
\begin{equation}
  E_{\rm opt}^{\rm num}
  \le
  \lVert\Phi_f\rVert_{\ell^2\to X} L_c(\mathcal N) \sqrt{\frac{2\varepsilon_{\rm tol}}{\mu}},
  \label{eq:algorithmic-state-error-bound}
\end{equation}
where the stopping tolerance bounds the objective gap under the assumptions above.
\end{remark}

\section{Numerical Experiments}
\label{sec:experiments}
\label{sec:numerical-experiments}

In this section, we present a series of numerical experiments to systematically assess the performance, accuracy, and robustness of the proposed multiscale kernel framework for joint state and parameter recovery. In alignment with the theoretical error decomposition established in Section~\ref{sec:theory} (specifically Proposition~\ref{prop:three-level-state-error-decomposition}, Theorem~\ref{thm:implementation-matched-infimal-multiscale}, and Proposition~\ref{prop:finite-penalty-relaxation}), our empirical investigation provides an empirical performance evaluation of how the operator-adapted multiscale kernel isolates fine-scale heterogeneities and how localized weak Galerkin formulations compare against discrete strong-form collocation.

To isolate the effects of the trial space and the PDE residual formulation, we compare four methods formed by pairing two trial spaces---the GMsFEM multiscale space (MS) and a smooth Gaussian RKHS (G)---with two residual formulations---the weak Galerkin residual (W) and strong-form point collocation (S)---across the same 50 synthetic test problems. Below, \cref{subsec:experimental-protocol} details the boundary value prototype, the synthetic coefficient parameterization, and the observation model. \cref{subsec:representation-information} formulates the comparative baselines and implementation parameters, and \cref{subsec:recovery-results} provides a comprehensive statistical and spatial error analysis.

\subsection{Experimental Setup and Problem Configuration}
\label{subsec:experimental-protocol}

We consider the steady-state scalar elliptic boundary value problem on the unit square $D = (0,1)^2$:
\begin{equation}
  -\nabla \cdot \left(\kappa^{\rm phys}(x)\nabla u^{\rm phys}(x)\right) = f(x) \quad \text{in } D,
  \qquad
  u^{\rm phys}(x) = 0 \quad \text{on } \partial D,
  \label{eq:num-forward-bvp}
\end{equation}
driven by a constant deterministic source term $f(x) \equiv 10$. The physical domain is discretized using two nested, conforming spatial meshes: a coarse triangulation $\mathcal{T}_H$ with mesh size $H = 1/10$, and a shape-regular fine triangulation $\mathcal{T}_h$ with mesh size $h = 1/100$, forming a Cartesian background partition of $100 \times 100$ square cells. The background coefficient geometry is defined by a reference permeability field $\kappa_{\rm ref} \in L^\infty(D)$ featuring high-conductivity channels and localized barrier inclusions embedded within a low-permeability matrix, yielding a severe contrast ratio of $\kappa_{\max}/\kappa_{\min} = 10^4$.

\paragraph{Admissible Coefficient Family.}
To explore parameter estimation within a controlled parameter domain, we specialize $\theta = \nu \in \mathbb{R}^{N_\theta}$ with $N_\theta = 20$ using the positive masked random Fourier feature (MRFF) family \citep{rahimi2007random}:
\begin{equation}
  \kappa_\nu(x)
  =
  \kappa_{\rm ref}(x)
  \exp\!\left(
    M(x)\sum_{r=1}^{N_\theta}
    \nu_r\cos(\boldsymbol\xi_r^\top x+b_r)
  \right),
  \qquad M(x) \equiv 1.
  \label{eq:mrff-coefficient}
\end{equation}
The spatial frequencies $\boldsymbol\xi_r \sim \mathcal{N}(0, I_2)$ and phase shifts $b_r \sim \operatorname{Unif}(0,2\pi)$ are sampled offline once and frozen across all samples. Consequently, the spatial basis, mask $M(x)$, and reference field $\kappa_{\rm ref}(x)$ remain invariant prior to online inversion. This setup focuses specifically on identifying the amplitude modulation vector $\nu \in \mathbb{R}^{20}$ on a compact parameter sublevel set (ensuring uniform ellipticity as in Assumption~A2), testing the recovery of multiplicative modulation over a prescribed high-contrast reference coefficient geometry rather than the discovery of unknown interface locations.

\paragraph{Data Synthesis and Inversion Protocol.}
We synthesize $N_{\rm samp}=50$ independent benchmark instances. For each realization $s \in \{1,\dots,N_{\rm samp}\}$, a ground-truth parameter vector is independently drawn from a centered Gaussian prior:
\begin{equation}
  \nu_s^\dagger \sim \mathcal{N}(0, \sigma_\nu^2 I_{20}), \qquad \sigma_\nu = 0.2.
\end{equation}
The target coefficient $\kappa_s^\dagger := \kappa_{\nu_s^\dagger}$ is evaluated at the barycenters of the fine cells and treated as cellwise constant on $\mathcal{T}_h$. The fine-scale ground-truth state $u_s^\dagger \in V_{h,0}$ is computed via the continuous conforming finite element Galerkin discretization of \cref{eq:num-forward-bvp} on $\mathcal{T}_h$ using piecewise linear basis functions.

Physical diagnostic measurements are extracted by deploying an array of $N_{\rm obs} = 100$ point sensors located strictly in the interior of $D$. The synthetic observation vector $y_s \in \mathbb{R}^{N_{\rm obs}}$ is perturbed by additive Gaussian white noise:
\begin{equation}
  y_s = \mathcal{O}(u_s^\dagger) + \eta_s, \qquad \eta_s \sim \mathcal{N}(0, \sigma_{\eta,s}^2 I_{N_{\rm obs}}),
  \label{eq:num-observation-data}
\end{equation}
where the noise standard deviation is scaled to one percent of the empirical target signal: $\sigma_{\eta,s} = 0.01 \times \operatorname{std}(\mathcal{O}(u_s^\dagger))$. All comparative methods share the identical dataset $\{(\kappa_s^\dagger, u_s^\dagger, y_s)\}_{s=1}^{N_{\rm samp}}$ under a matched computational budget. Discretization settings and inversion parameters are summarized in \cref{tab:shared-experimental-protocol}.

\begin{table}[H]
  \centering
  \caption{Benchmark discretization, observational dimensions, and optimization parameters.}
  \label{tab:shared-experimental-protocol}
  \footnotesize
  \renewcommand{\arraystretch}{1.15}
  \begin{tabular}{@{}>{\raggedright\arraybackslash}p{0.24\textwidth}
                      p{0.72\textwidth}@{}}
    \toprule
    Category & Implementation Details \\
    \midrule
    Forward Problem \& Domain
      & Unit square $D=(0,1)^2$, homogeneous Dirichlet boundary conditions ($u=0$ on $\partial D$), source term $f(x) \equiv 10$. Reference coefficient $\kappa_{\rm ref}$ with contrast ratio $\kappa_{\max}/\kappa_{\min} = 10^4$. Background grid of $100\times100$ fine cells; coarse GMsFEM partition of $10\times10$ elements. \\
    Observation Setup \& Noise
      & $N_{\rm samp}=50$ trials; $N_{\rm obs}=100$ strict-interior observation sensors. Additive Gaussian noise with $\sigma_{\eta,s} = 0.01 \times \operatorname{std}(\mathcal{O}(u_s^\dagger))$. Problem dimensions: parameter $N_\theta=20$, free-state $N_c=324$, and physics constraints $N_{\rm weak}=324$. \\
    Optimization Protocol
      & State and parameter regularizers fixed at $\alpha=\lambda_u=10^{-2}$. Physics weighting scales: $\beta=10^{-4}$ (multiscale/weak), $\beta=4.110\times10^{-4}$ (Gaussian/weak), and $\beta=9.715\times10^{-2}$ (both discrete strong-form cells). Variable projection initialized at $\nu^0=0$, solved via unconstrained L-BFGS-B with a fixed termination budget of $K_{\max}=10$ iterations. \\
    \bottomrule
  \end{tabular}
\end{table}

\paragraph{Evaluation Metrics.}
Upon the optimizer returning the terminal parameter estimate $\hat{\nu}_s$ and the associated projected multiscale state coordinates $\hat{c}_s = c_*(\hat{\nu}_s)$, the physical reconstruction fields are evaluated as $u_{{\rm rec},s} = \Phi_f \hat{c}_s$ and $\kappa_{{\rm rec},s} = \kappa_{\hat{\nu}_s}$. We quantify recovery accuracy via relative discrete $L^2$ error metrics:
\begin{equation}
  e_u^{(s)}
  :=
  \frac{\lVert u_{{\rm rec},s}-u_s^\dagger\rVert_{2,h}}
       {\lVert u_s^\dagger\rVert_{2,h}},
  \qquad
  e_{\log\kappa}^{(s)}
  :=
  \frac{\lVert\log\kappa_{{\rm rec},s}-\log\kappa_s^\dagger\rVert_{2,\mathcal{T}_h}}
       {\lVert\log\kappa_s^\dagger\rVert_{2,\mathcal{T}_h}}.
  \label{eq:experimental-relative-errors}
\end{equation}
Here, $\lVert\cdot\rVert_{2,h}$ denotes the Euclidean vector norm evaluated across all interior fine-grid vertices $\mathcal{N}_h^0$, whereas $\lVert\cdot\rVert_{2,\mathcal{T}_h}$ denotes the cellwise piecewise constant $L^2(D)$ norm over the fine background grid. In the subsequent comparative tables, we report the sample mean and empirical standard deviation over the $N_{\rm samp}=50$ trials.

\subsection{Trial Spaces and PDE Residual Formulations}
\label{subsec:representation-information}

The comparative framework crosses two distinct trial spaces, indexed by $X\in\{\mathrm{MS},\mathrm G\}$, with two PDE residual formulations, indexed by $Y\in\{\mathrm W,\mathrm S\}$, as summarized in \cref{tab:comparison-factors}. Here, $(X,Y)=(\mathrm{MS},\mathrm W)$ designates the proposed multiscale kernel framework.

\begin{table}[H]
  \centering
  \caption{Methodological factors in the four-way comparison.}
  \label{tab:comparison-factors}
  \small
  \renewcommand{\arraystretch}{1.15}
  \begin{tabular}{@{}llp{0.54\textwidth}@{}}
    \toprule
    Factor Category & Identifier & Mathematical Formulation \& Role \\
    \midrule
    Trial space ($X$) & $\mathrm{MS}$
      & Operator-adapted GMsFEM multiscale basis $\Phi_{\rm ms}$ encoding $\kappa_{\rm ref}$ \\
      & $\mathrm G$
      & Standard translation-invariant Gaussian kernel basis $\Phi_{\rm G}$ \\
    \addlinespace[0.5ex]
    PDE residual formulation ($Y$) & $\mathrm W$
      & Localized weak Galerkin variational balance via test functions \\
      & $\mathrm S$
      & Discrete strong-form collocation residual at interior Halton points \\
    \bottomrule
  \end{tabular}
\end{table}

\paragraph{Trial Spaces ($X$).}
Both candidate spaces are configured with an identical dimension of $N_c = 324$ free degrees of freedom, parameterizing candidate states as $u_X(x) = \Phi_X(x) c$ with coordinate vector $c \in \mathbb{R}^{324}$.
\begin{enumerate}[leftmargin=1.5em, label=(\roman*)]
  \item \textbf{GMsFEM-induced Multiscale Space ($X = \mathrm{MS}$):}
  We set $\Phi_{\rm MS} = \Phi_f$, where $\Phi_f(x) = [\Phi_{f,1}(x),\dots,\Phi_{f,324}(x)] \in \mathbb{R}^{1\times 324}$ comprises the interior conforming GMsFEM basis functions constructed offline in \cref{subsec:multiscale-state-kernel} from the reference coefficient $\kappa_{\rm ref}$. This basis induces the operator-adapted reproducing kernel:
  \begin{equation}
    K_{\rm ms}(x,z) = \sum_{j=1}^{324} \Phi_{f,j}(x)\Phi_{f,j}(z) = \Phi_f(x)\Phi_f(z)^\top, \qquad x,z \in D.
    \label{eq:comp-ms-kernel}
  \end{equation}
  By construction, $\mathcal{H}_{K_{\rm ms}} = V_{{\rm ms},0} \subset H_0^1(D)$, naturally enforcing zero trace on $\partial D$ while capturing rough interfacial gradients without scale separation.

  \item \textbf{Matched Gaussian Space ($X = \mathrm{G}$):}
  As a smooth baseline, we consider the isotropic squared-exponential Gaussian kernel:
  \begin{equation}
    K_{\rm G}(x,z) = \exp\!\left(-\frac{\lVert x-z\rVert_2^2}{2\sigma_{\rm G}^2}\right), \qquad \sigma_{\rm G} = 0.2.
    \label{eq:comp-gaussian-kernel}
  \end{equation}
  To construct a conforming state representation of matched dimension, we deploy an array of $18 \times 18 = 324$ distinct collocation centers $\{z_j\}_{j=1}^{324}$ on a uniform lattice strictly within the interior of $D$. The basis functions $\Phi_{\mathrm G} = [\phi_1^{\mathrm G},\dots,\phi_{324}^{\mathrm G}]$ are obtained by multiplying each kernel section by a smooth boundary-distance mollifier ensuring zero trace on $\partial D$. Furthermore, $\Phi_{\mathrm G}$ is normalized such that its average kernel diagonal over all interior fine-grid nodes matches that of the multiscale kernel, i.e., $\sum_{p\in\mathcal{N}_h^0} K_{\rm G}(p,p) = \sum_{p\in\mathcal{N}_h^0} K_{\rm ms}(p,p)$.
\end{enumerate}

\paragraph{PDE Residual Formulations ($Y$).}
The physical balance governed by $\kappa_\nu$ is enforced through either a localized integral weak formulation or a discrete strong-form differential residual.
\begin{enumerate}[leftmargin=1.5em, label=(\roman*)]
  \item \textbf{Weak Galerkin Formulation ($Y = \mathrm{W}$):}
  At the continuous functional level, the weak residual functional evaluates variational compatibility:
  \begin{equation}
    \mathcal{R}_{\mathrm W}(u_X,\nu; v)
    := \int_D \kappa_\nu(x)\nabla u_X(x) \cdot \nabla v(x)\,dx - \langle f, v \rangle_{H^{-1}, H_0^1},
    \qquad v \in H_0^1(D).
    \label{eq:comp-weak-residual-cont}
  \end{equation}
  Choosing the 324 basis functions of the respective space as test functions ($\{v_i = \Phi_{X,i}\}_{i=1}^{324}$) and approximating integrals via fine-grid numerical quadrature, we obtain the algebraic residual $r_{X,\mathrm W}(c,\nu) := A_X(\nu)c - b_X \in \mathbb{R}^{324}$, where $[A_X(\nu)]_{ij} \approx a_{\kappa_\nu}(\Phi_{X,j}, \Phi_{X,i})$ and $[b_X]_i \approx \ell(\Phi_{X,i})$. This formulation requires only $H^1$-regularity, enforcing discrete variational balance without evaluating distributional derivatives pointwise across coefficient jumps.

  \item \textbf{Discrete Strong-form Collocation ($Y = \mathrm{S}$):}
  Collocation approaches evaluate the strong-form differential residual:
  \begin{equation}
    \mathcal{R}_{\mathrm S}(u_X,\nu)(x) := -\nabla \cdot \left(\kappa_\nu(x)\nabla u_X(x)\right) - f(x).
    \label{eq:comp-strong-residual-cont}
  \end{equation}
  Because $\kappa_\nu$ is high-contrast and discontinuous, classical strong derivatives do not exist pointwise across interfaces. To construct a concrete discrete baseline on $\mathcal{T}_h$, let $A_h(\kappa_\nu) \in \mathbb{R}^{N_h \times N_h}$ denote the fine-scale stiffness matrix, $M_h \in \mathbb{R}^{N_h \times N_h}$ the lumped fine-grid mass matrix, and $P_X \in \mathbb{R}^{N_h \times 324}$ the prolongation matrix of the trial basis $\Phi_X$. We define the discrete collocation operator by $S_X(\nu) := E_{\rm col} M_h^{-1} A_h(\kappa_\nu) P_X \in \mathbb{R}^{324 \times 324}$, where $E_{\rm col} \in \mathbb{R}^{324 \times N_h}$ evaluates or interpolates fine-grid nodal balances onto 324 quasi-random interior collocation points $\{x_k^{\rm col}\}_{k=1}^{324}$ generated via a deterministic two-dimensional Halton sequence. The corresponding discrete load vector is $f_{\mathrm S} := E_{\rm col} M_h^{-1} b_h \in \mathbb{R}^{324}$, yielding the algebraic collocation residual $r_{X,\mathrm S}(c,\nu) := S_X(\nu)c - f_{\mathrm S} \in \mathbb{R}^{324}$.
\end{enumerate}

\paragraph{Unified Objective Functional and Variable Projection.}
For each instance $s \in \{1,\dots,N_{\rm samp}\}$ and methodological cell $(X,Y)$, joint inversion minimizes the regularized functional:
\begin{equation}
  J_{X,Y}^{(s)}(c,\nu)
  :=
  \frac{1}{2\gamma_s^2}
  \bigl\lVert H_X^{(s)}c - y_s\bigr\rVert_2^2
  +
  \frac{1}{2\beta_{X,Y}^2}
  \bigl\lVert r_{X,Y}(c,\nu)\bigr\rVert_2^2
  +
  \frac{\lambda_u}{2}\lVert c\rVert_2^2
  +
  \frac{\alpha}{2}\lVert\nu\rVert_2^2,
  \label{eq:experimental-common-comparator-objective}
\end{equation}
where $H_X^{(s)} \in \mathbb{R}^{100 \times 324}$ denotes the state observation matrix, $\gamma_s = \sigma_{\eta,s}$, and $\lambda_u = \alpha = 10^{-2}$. The penalty weights $\beta_{X,Y}$ are fixed a priori to balance the dimensional scale of the physics residual against the load vector, as specified in \cref{tab:shared-experimental-protocol}.

For each benchmark instance $s$ and methodological cell $(X,Y)$, the conditionally linear state coordinates are eliminated via the thin QR factorization of the augmented design matrix following \cref{subsec:variable-projection}. The reduced objective $\varphi_{X,Y}^{(s)}( \nu)$ is minimized using unconstrained L-BFGS-B initialized at $\nu^0=0$, with its gradient computed via the adjoint-free formula \cref{eq:reduced-objective-gradient} under a fixed budget of $K_{\max}=10$ iterations.

\subsection{Recovery Results and Comparative Analysis}
\label{subsec:recovery-results}

\cref{tab:factor-ablation} compiles the sample means and standard deviations of the relative state error $e_u$ and log-coefficient error $e_{\log\kappa}$ across the $N_{\rm samp}=50$ test instances under the unified optimization protocol.

\begin{table}[H]
  \centering
  \caption{Statistical performance across the four-way comparison ($N_{\rm samp}=50$).} Entries report the sample mean (with sample standard deviation in parentheses). The lowest mean in each error column is highlighted in bold.
  \label{tab:factor-ablation}
  \footnotesize
  \renewcommand{\arraystretch}{1.15}
  \setlength{\tabcolsep}{5pt}
  \begin{tabular}{@{}>{\raggedright\arraybackslash}p{0.25\textwidth}
                      >{\raggedright\arraybackslash}p{0.22\textwidth}
                      rr@{}}
    \toprule
    Trial Space ($X$) & PDE Residual Formulation ($Y$)
      & State Error ($e_u$) & Log-Coefficient Error ($e_{\log\kappa}$) \\
    \midrule
    GMsFEM-induced multiscale
      & Weak Galerkin (proposed)
      & \textbf{0.0237} (0.0076)
      & \textbf{0.0443} (0.0336) \\
    GMsFEM-induced multiscale
      & Discrete strong
      & 0.2695 (0.0822)
      & 0.1654 (0.0892) \\
    \addlinespace[0.4ex]
    Matched Gaussian
      & Weak Galerkin
      & 0.1635 (0.0287)
      & 1.5557 (0.1275) \\
    Matched Gaussian
      & Discrete strong
      & 0.1659 (0.0175)
      & 0.3035 (0.0433) \\
    \bottomrule
  \end{tabular}
\end{table}

\Cref{tab:factor-ablation} reports the mean and standard deviation of the errors for each configuration. All four configurations used the same 50 test instances and a maximum of ten L-BFGS-B iterations. With the weak residual, replacing the Gaussian trial space with the multiscale space reduced the mean relative state error from 16.35\% to 2.37\%. Within the multiscale space, replacing the discrete strong residual with the weak residual reduced that error from 26.95\% to 2.37\%. The corresponding mean log-coefficient error decreased from 16.54\% to 4.43\%. These comparisons use a fixed basis dimension and a single reference contrast.

\begin{figure}[H]
  \centering
  \begin{subfigure}[t]{0.45\linewidth}
    \centering
    \includegraphics[width=\linewidth]{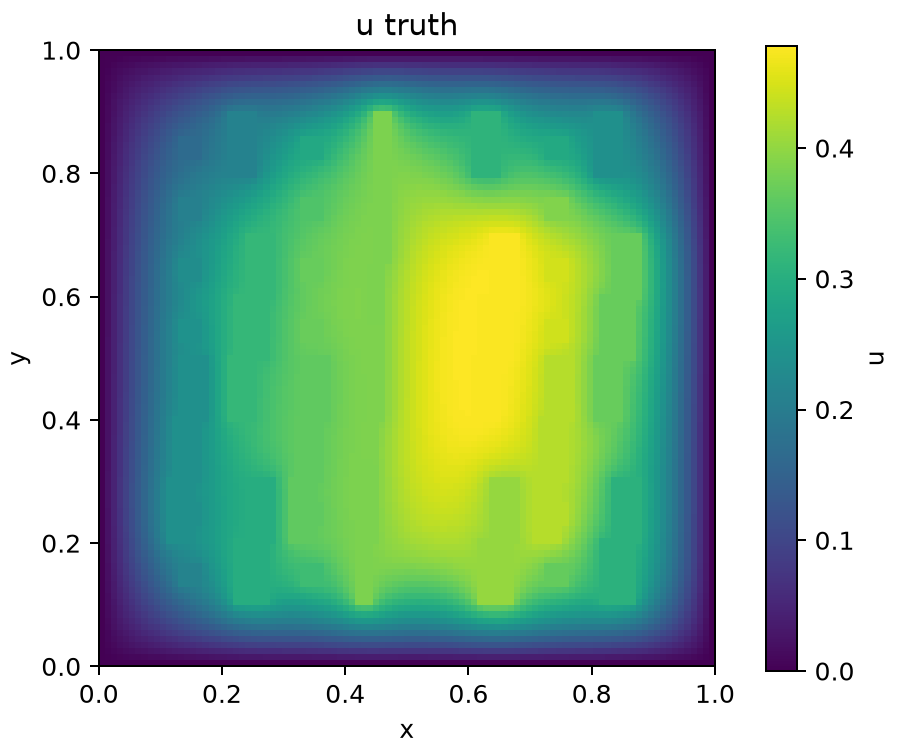}
    \caption{Fine-grid ground-truth state $u_0^\dagger$.}
  \end{subfigure}\hfill
  \begin{subfigure}[t]{0.45\linewidth}
    \centering
    \includegraphics[width=\linewidth]{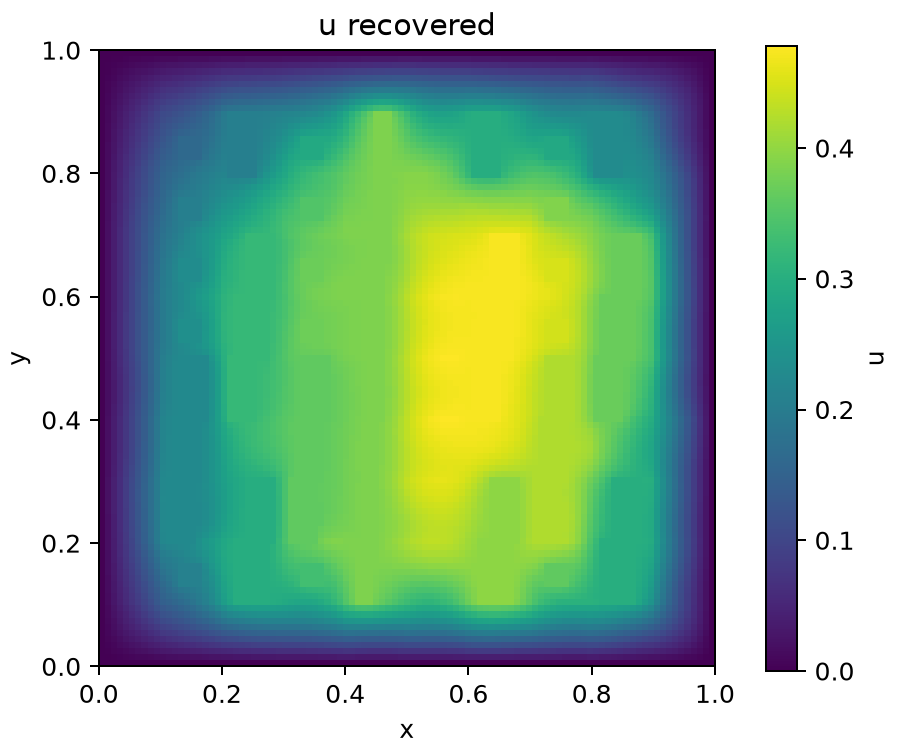}
    \caption{Recovered multiscale state $u_{\rm rec, 0}$.}
  \end{subfigure}

  \vspace{1ex}
  \begin{subfigure}[t]{0.45\linewidth}
    \centering
    \includegraphics[width=\linewidth]{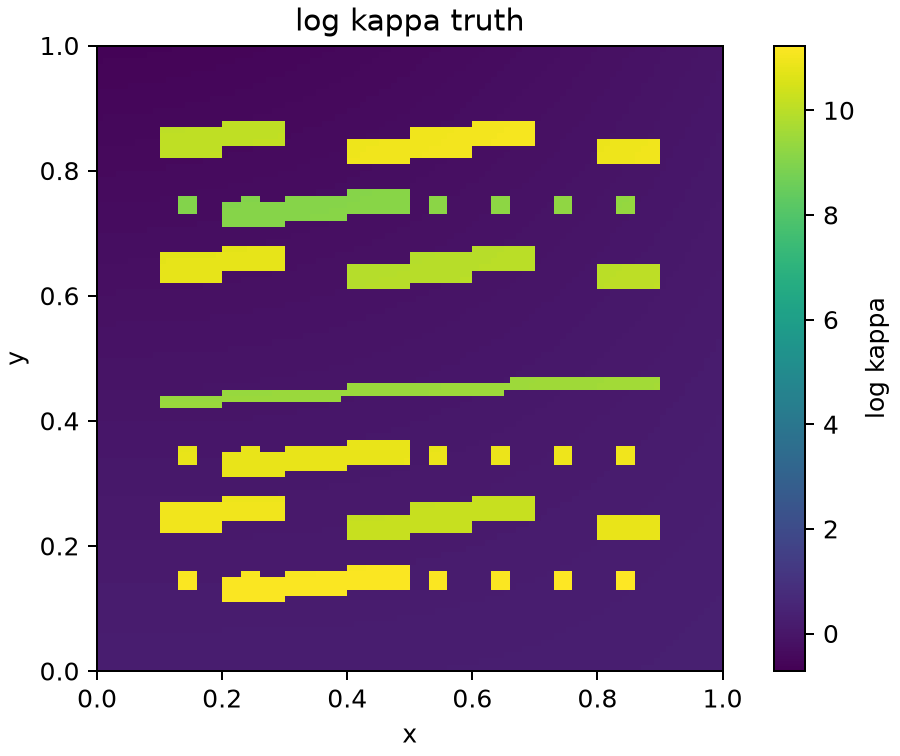}
    \caption{Target log-coefficient $\log\kappa_0^\dagger$.}
  \end{subfigure}\hfill
  \begin{subfigure}[t]{0.45\linewidth}
    \centering
    \includegraphics[width=\linewidth]{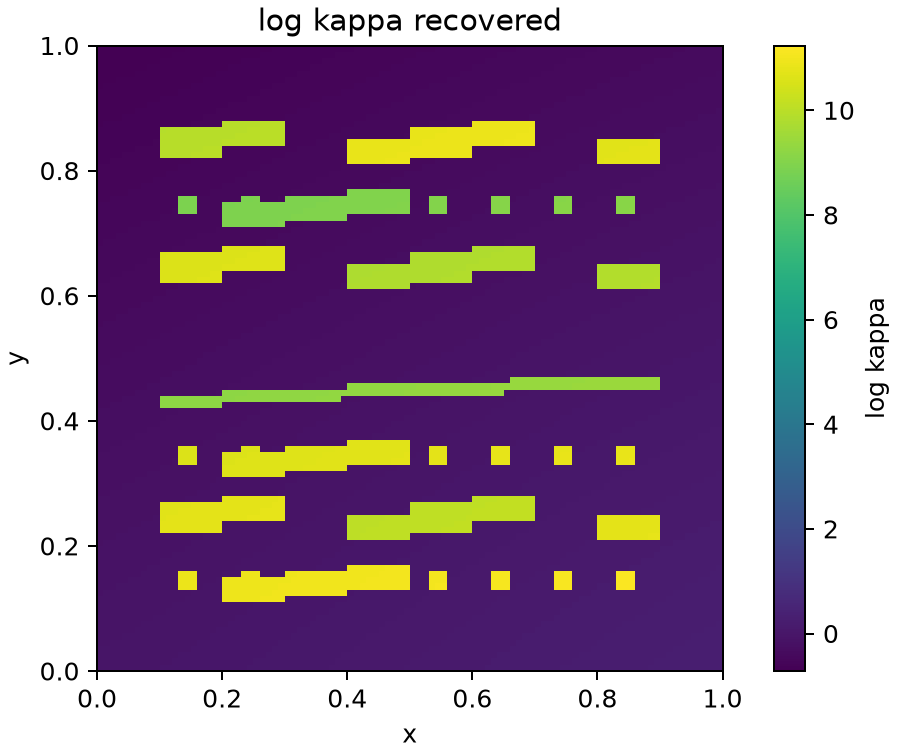}
    \caption{Recovered effective log-coefficient $\log\kappa_{\rm rec, 0}$.}
  \end{subfigure}
  
  \caption{Spatial reconstruction profiles for representative test sample $s=0$ using the proposed multiscale/weak framework (MS-W). The individual relative errors ($e_u^{(0)} = 0.0241$, $e_{\log\kappa}^{(0)} = 0.0452$) closely match the ensemble sample means reported in \cref{tab:factor-ablation}.}
  \label{fig:end-to-end-representative}
\end{figure}

\paragraph{Spatial Reconstruction and Error Localization.}
To complement the statistical summary, \cref{fig:end-to-end-representative} visualizes the spatial fields for sample $s=0$, selected because its individual errors ($e_u^{(0)} = 0.0241$, $e_{\log\kappa}^{(0)} = 0.0452$) lie near the ensemble average.

Panels (a) and (b) show the target and recovered states. Panels (c) and (d) show the target and recovered log-coefficients. The reconstructed coefficient varies in amplitude over the prescribed reference geometry, with visible differences near some inclusion boundaries.

\section{Conclusion}
\label{sec:conclusion}

We presented a method for estimating the state and a parameterized diffusion coefficient in a high-contrast elliptic equation. The method uses a GMsFEM basis computed from a prescribed reference coefficient and a weak Galerkin residual. Variable projection reduces the optimization to the coefficient parameters. Our error estimates distinguish multiscale approximation, penalty relaxation, and parameter estimation.

In the synthetic tests, the multiscale space with the weak residual gave the lowest mean state and log-coefficient errors among the four configurations. The tests recover coefficient modulation over prescribed channel and barrier locations. Recovery when the reference geometry differs from the target remains an open question for this formulation.

\section*{Acknowledgement}

Eric T. Chung's work is partially supported by the Hong Kong RGC General Research Fund (Project Nos. 14304525 and 14305624).

\section*{Data Availability Statement}
The data and code supporting the numerical findings of this study are available from the corresponding author upon reasonable request.

\section*{Conflict of Interest}
The authors declare that they have no known competing financial interests or personal relationships that could have appeared to influence the work reported in this paper.

\appendix


\section{Finite-Dimensional Recovery, MAP, and Variable Projection}
\label{app:finite-dimensional-recovery}

This appendix provides the proofs for the structural and algebraic results established in \cref{subsec:finite-rank-projected-state}.

\subsection{Finite-Rank State Space and Optional Joint MAP}
\label{app:finite-rank-joint-map-proof}

\begin{proof}[Proof of Proposition~\ref{prop:finite-rank-joint-map}]
We first verify the deterministic reproducing kernel properties. By assumption
\ref{ass:features}, the basis functions $\{\Phi_{f,j}\}_{j=1}^{N_c}$ form a linearly
independent basis of $V_{\rm ms,0}$. Thus, the linear synthesis operator
$T:\mathbb R^{N_c}\to V_{\rm ms,0}$, defined by $T c := \Phi_f c = \sum_{j=1}^{N_c} c_j \Phi_{f,j}$,
is a linear bijection. Endow $V_{\rm ms,0}$ with the inner product
$\langle Tc, Td\rangle_{K_{\rm ms}} := c^\top d$ for all $c, d\in\mathbb R^{N_c}$. The induced norm satisfies
$\lVert u_c\rVert_{K_{\rm ms}}^2 = c^\top c = \lVert c\rVert_2^2$, establishing an isometric isomorphism
between $(\mathbb R^{N_c}, \lVert\cdot\rVert_2)$ and $(V_{\rm ms,0}, \lVert\cdot\rVert_{K_{\rm ms}})$.

Recalling the definition of the multiscale kernel $K_{\rm ms}(x,x') = \Phi_f(x)\Phi_f(x')^\top$,
we observe that for each fixed $x'\in D$, the section $K_{\rm ms}(\cdot, x') = \Phi_f(\cdot)\Phi_f(x')^\top$
belongs to $V_{\rm ms,0}$ with coordinate vector $d_{x'} = \Phi_f(x')^\top\in\mathbb R^{N_c}$. For any
$u_c = \Phi_f c\in V_{\rm ms,0}$, the reproducing property holds:
\[
  \langle u_c, K_{\rm ms}(\cdot, x')\rangle_{K_{\rm ms}}
  =
  \langle T c, T d_{x'}\rangle_{K_{\rm ms}}
  =
  c^\top d_{x'}
  =
  \Phi_f(x') c
  =
  u_c(x'),
  \qquad\forall x'\in D.
\]
By the uniqueness of reproducing kernels, $\mathcal H_{K_{\rm ms}} = V_{\rm ms,0}$ and its native norm
is $\lVert u_c\rVert_{K_{\rm ms}} = \lVert c\rVert_2$, which proves \cref{eq:finite-rank-rkhs-identity}.
The identity $\frac{\lambda_u}{2}\lVert c\rVert_2^2 = \frac{\lambda_u}{2}\lVert u_c\rVert_{K_{\rm ms}}^2$
follows immediately.

To prove the statistical MAP equivalence, consider the conditionally independent Gaussian
likelihood models in \cref{eq:observation-likelihood}:
\[
  \pi(y\mid c) \propto \exp\left(-\frac{1}{2}\lVert H_f c - y\rVert_{W_y}^2\right),
  \qquad
  \pi(b_f\mid c,\theta) \propto \exp\left(-\frac{1}{2}\lVert A_f(\theta)c - b_f\rVert_{W_W}^2\right).
\]
The independent prior densities on $(c,\theta)$ defined in \cref{eq:coordinate-priors} satisfy
\[
  \pi_0(c) \propto \exp\left(-\frac{\lambda_u}{2}\lVert c\rVert_2^2\right),
  \qquad
  \pi_0(\theta) \propto \exp\left(-R_\theta(\theta)\right).
\]
Let $\mu_\Theta$ denote a Radon reference measure supported on the closed parameter domain $\Theta$ (coinciding with the standard Lebesgue measure whenever $\Theta$ has non-empty interior in $\mathbb R^{N_\theta}$).
Rather than relying on coercivity alone, we explicitly assume that the prior density is integrable, satisfying $0 < Z_\theta := \int_\Theta \exp(-R_\theta(\theta))\,d\mu_\Theta(\theta) < \infty$, which holds unconditionally for standard quadratic Tikhonov regularization $R_\theta(\theta) = \frac{\alpha}{2}\lVert\theta\rVert_2^2$ on $\mathbb R^{N_\theta}$.
Thus, $\pi_0(\theta) = Z_\theta^{-1}\exp(-R_\theta(\theta))$ defines a proper prior probability density on $\Theta$. By Bayes' theorem, the joint
posterior distribution with respect to the product reference measure $dc \otimes d\mu_\Theta(\theta)$ on $\mathbb R^{N_c}\times\Theta$
satisfies:
\[
  \pi(c,\theta\mid y, b_f)
  \propto
  \pi(y\mid c)\,\pi(b_f\mid c,\theta)\,\pi_0(c)\,\pi_0(\theta)
  =
  \exp\left(-J(c,\theta)\right),
\]
where $J(c,\theta)$ is the variational functional defined in \eqref{eq:joint-regularized-objective}.
The joint likelihood is bounded above by $1$, and the joint prior is integrable; hence the
posterior measure is proper. Taking the negative logarithm reveals that any global
minimizer $(c_\star, \theta_\star)$ of $J(c,\theta)$ maximizes the posterior density, coinciding
with a joint MAP estimate. Because the coupling term $A_f(\theta)c$ is linear in the state coordinates $c$ but generally nonlinear in the parameter coordinates $\theta$ (as in exponential or non-affine coefficient parameterizations), and $R_\theta$ may be non-quadratic, the joint posterior is in general non-Gaussian.
\end{proof}

\subsection{Projected State and Reduced Objective}
\label{app:projected-state-reduced-objective-proof}

\begin{proof}[Proof of Proposition~\ref{prop:projected-state-reduced-objective}]
For any fixed parameter $\theta\in\Theta$, the quadratic state subproblem functional is
$J(c,\theta) = \frac{1}{2}\lVert D(\theta)c - d\rVert_2^2 + R_\theta(\theta)$. The Hessian with respect to $c$
is given by the augmented Gram matrix:
\[
  \nabla_c^2 J(c,\theta)
  =
  D(\theta)^\top D(\theta)
  =
  H_f^\top W_y H_f + A_f(\theta)^\top W_W A_f(\theta) + \lambda_u I_{N_c}.
\]
Because $W_y\succ0$ and $W_W\succ0$, the terms $H_f^\top W_y H_f$ and $A_f(\theta)^\top W_W A_f(\theta)$
are positive semi-definite. Therefore,
\begin{equation}
  z^\top D(\theta)^\top D(\theta) z
  =
  \lVert H_f z\rVert_{W_y}^2 + \lVert A_f(\theta) z\rVert_{W_W}^2 + \lambda_u\lVert z\rVert_2^2
  \ge
  \lambda_u\lVert z\rVert_2^2,
  \qquad\forall z\in\mathbb R^{N_c}.
  \label{eq:coercivity-gram}
\end{equation}
This establishes \cref{eq:design-matrix-spectral-bound}, showing that $D(\theta)$ has full column
rank $N_c$ and that $c\mapsto J(c,\theta)$ is strictly and $\lambda_u$-strongly convex. Hence,
$J(c,\theta)$ admits a unique global minimizer $c_\star(\theta)$. For the thin QR factorization
$D(\theta)=Q_D(\theta)R_D(\theta)$, the upper-triangular factor $R_D(\theta)\in\mathbb R^{N_c\times N_c}$
is nonsingular, and the solution to the triangular system $R_D(\theta)c_\star(\theta)=Q_D(\theta)^\top d$
coincides with the unique state.

The first-order optimality condition $\nabla_c J(c_\star(\theta),\theta)=0$ yields the normal equations:
\begin{equation}
  M(\theta) c_\star(\theta) = g(\theta),
  \label{eq:normal-equations-rep}
\end{equation}
where
\[
  M(\theta) := D(\theta)^\top D(\theta) \succeq \lambda_u I_{N_c},
  \qquad
  g(\theta) := H_f^\top W_y y + A_f(\theta)^\top W_W b_f.
\]
By assumption \ref{ass:parameter-regularity}, $\theta\mapsto A_f(\theta)$ is $C^1$. Consequently,
$M(\theta)$ and $g(\theta)$ are $C^1$ mappings on an open neighborhood of $\Theta$. Because the
matrix inversion operator $M\mapsto M^{-1}$ is smooth on the cone of symmetric positive definite
matrices, the solution mapping $c_\star(\theta) = M(\theta)^{-1}g(\theta)$ is $C^1$. It follows
immediately that the reduced objective $\varphi(\theta) = J(c_\star(\theta),\theta)$ is $C^1$ on $\Theta$.

To establish existence of a global minimizer for $\varphi$, we observe that all least-squares
mismatch terms are non-negative, which gives the lower bound:
\[
  \varphi(\theta) = \frac{1}{2}\lVert D(\theta)c_\star(\theta)-d\rVert_2^2 + R_\theta(\theta) \ge R_\theta(\theta).
\]
By coercivity of $R_\theta$ on $\Theta$ (assumption \ref{ass:parameter-regularity}),
$\varphi(\theta)\to\infty$ as $\|\theta\|_2\to\infty$. Let $\alpha := \inf_{\theta\in\Theta}\varphi(\theta)$. For any
sublevel set $S_\gamma := \{\theta\in\Theta : \varphi(\theta)\le\gamma\}$ with $\gamma > \alpha$, coercivity implies that
$S_\gamma$ is bounded. Because $\Theta$ is closed and $\varphi$ is continuous, $S_\gamma$ is compact and
nonempty. By the Weierstrass extreme value theorem, $\varphi$ achieves its infimum on $S_\gamma$,
guaranteeing the existence of at least one global minimizer $\theta_\star\in\Theta$.

We now derive the explicit reduced gradient formula \cref{eq:reduced-objective-gradient}.
By the multivariate chain rule, the gradient of $\varphi(\theta) = J(c_\star(\theta),\theta)$ is:
\begin{equation}
  \nabla_\theta \varphi(\theta)
  =
  \nabla_\theta J(c,\theta)\Big|_{c=c_\star(\theta)}
  +
  \left(D_\theta c_\star(\theta)\right)^\top \nabla_c J(c,\theta)\Big|_{c=c_\star(\theta)}.
  \label{eq:chain-rule-expansion}
\end{equation}
Because $c_\star(\theta)$ uniquely minimizes $J(\cdot,\theta)$, the first-order stationarity condition $\nabla_c J(c_\star(\theta),\theta) = 0$ holds identically. The second term in \cref{eq:chain-rule-expansion} therefore vanishes, leaving:
\[
  \nabla_\theta \varphi(\theta) = \nabla_\theta J(c_\star(\theta),\theta).
\]
Differentiating $J(c,\theta) = \frac{1}{2}\lVert H_f c - y\rVert_{W_y}^2 + \frac{1}{2}\lVert A_f(\theta)c - b_f\rVert_{W_W}^2 + \frac{\lambda_u}{2}\lVert c\rVert_2^2 + R_\theta(\theta)$
partially with respect to $\theta_r$ at fixed $c = c_\star(\theta)$ yields:
\[
  \frac{\partial\varphi(\theta)}{\partial\theta_r}
  =
  \left(A_f(\theta)c_\star(\theta) - b_f\right)^\top W_W \frac{\partial A_f(\theta)}{\partial\theta_r} c_\star(\theta)
  +
  \frac{\partial R_\theta(\theta)}{\partial\theta_r},
\]
which is precisely \cref{eq:reduced-objective-gradient} since $r_W(c_\star(\theta),\theta) = A_f(\theta)c_\star(\theta) - b_f$.

Finally, if $\theta\mapsto A_f(\theta)$ and $R_\theta(\theta)$ are real analytic on a neighborhood of $\Theta$,
then $M(\theta)$ and $g(\theta)$ are real analytic. Since matrix inversion is analytic on the SPD
cone, $c_\star(\theta) = M(\theta)^{-1}g(\theta)$ is real analytic, which implies that $\varphi(\theta)$
is a composition of analytic functions and is therefore real analytic.
\end{proof}


\section{Multiscale Kernel Approximation Bounds}
\label{app:multiscale-approximation}

This appendix provides the complete proof of \cref{thm:implementation-matched-infimal-multiscale},
analyzing the spectral tail decay, the boundary-free admissible decompositions, and the contrast-independent
stability of the multiscale nodal assembly.

\subsection{Proof of Theorem~\ref{thm:implementation-matched-infimal-multiscale}}
\label{app:multiscale-approximation-proof}

The proof proceeds in three main steps: (i) establishing local spectral tail estimates; (ii) verifying the
analytical properties of the multiscale error functional under the admissible decomposition hypothesis;
and (iii) proving the uniform bound on the stability constant $C_{\times,i}^{\rm ref}$ of the partition-of-unity nodal interpolation operator and assembling the global estimate.

\subsubsection{Step 1: Local spectral tails}
For each coarse patch $\omega_i$ ($i\in\mathcal I_0$), the snapshot space $W_i = V_{\rm snap}^{(i)}$ is
finite-dimensional, and $s_i^+$ is symmetric positive definite because the integrand contains $1>0$.
Thus, there exists an $s_i^+$-orthonormal basis of generalized eigenfunctions $\{\psi_{i,\ell}\}_{\ell=1}^{\dim W_i}$
satisfying
\begin{equation}
  a_i^{\rm ref}(\psi_{i,\ell}, z) = \lambda_{i,\ell} s_i^+(\psi_{i,\ell}, z)
  \qquad\forall z\in W_i,
  \label{eq:app-eigen-rel}
\end{equation}
with eigenvalues $0 \le \lambda_{i,1} \le \lambda_{i,2} \le \cdots$. On a connected patch, the kernel of
$a_i^{\rm ref}$ consists solely of constant functions. Because $L_i \ge 1$, the first omitted eigenvalue satisfies
$\lambda_{i,L_i+1} > 0$.

Any local state $z\in W_i$ expands uniquely as $z = \sum_{\ell} c_{i,\ell} \psi_{i,\ell}$ with coordinates
$c_{i,\ell} = s_i^+(z, \psi_{i,\ell})$. The truncated component $e_i := z - \Pi_i^{(L_i)} z = \sum_{\ell \ge L_i+1} c_{i,\ell} \psi_{i,\ell}$
satisfies the spectral tail relations:
\begin{equation}
  s_i^+(e_i, e_i)
  =
  \sum_{\ell \ge L_i+1} c_{i,\ell}^2
  \le
  \frac{1}{\lambda_{i,L_i+1}} \sum_{\ell \ge L_i+1} \lambda_{i,\ell} c_{i,\ell}^2
  \le
  \frac{1}{\lambda_{i,L_i+1}} a_i^{\rm ref}(z, z).
  \label{eq:appendix-local-mass-tail}
\end{equation}
Furthermore, applying $\mathcal A_i \psi_{i,\ell} = \lambda_{i,\ell} \psi_{i,\ell}$ and $s_i^+$-orthonormality yields:
\begin{equation}
  a_i^{\rm ref}(e_i, e_i)
  =
  \sum_{\ell \ge L_i+1} \lambda_{i,\ell} c_{i,\ell}^2
  \le
  \frac{1}{\lambda_{i,L_i+1}} \sum_{\ell \ge L_i+1} \lambda_{i,\ell}^2 c_{i,\ell}^2
  \le
  \frac{1}{\lambda_{i,L_i+1}} s_i^+(\mathcal A_i z, \mathcal A_i z).
  \label{eq:appendix-local-stiffness-tail}
\end{equation}
Summing \cref{eq:appendix-local-mass-tail} and \cref{eq:appendix-local-stiffness-tail} confirms the local
spectral bound:
\begin{equation}
  a_i^{\rm ref}(e_i, e_i) + s_i^+(e_i, e_i) \le \frac{1}{\lambda_{i,L_i+1}} \mathfrak R_i(z).
  \label{eq:app-tail-sum}
\end{equation}

\subsubsection{Step 2: Admissible decompositions and snapshot space properties}
Let $\mathcal I = \mathcal I_0 \cup \mathcal I_\partial$ (with $\mathcal I_0\cap \mathcal I_\partial=\emptyset$) be the index set of the coarse vertices in $D$ and on $\partial D$, respectively.
By construction, the full initial partition of unity $\{\chi_i^{\rm ms}\}_{i\in\mathcal I}$ satisfies $\sum_{i\in\mathcal I} \chi_i^{\rm ms}(p) = 1$ at every fine-grid node $p\in\overline D$.
Regarding the discrete stiffness operator $A_h$, we note that under cellwise constant diffusion $\kappa_{\rm ref}$, the off-diagonal entry corresponding to an interior fine edge $e = T_1 \cap T_2$ evaluates to $A_{h,ij} = -\frac{1}{2}(\kappa_1 \cot \alpha_1 + \kappa_2 \cot \alpha_2)$, where $\alpha_1, \alpha_2$ are the angles opposite to $e$.
Under severe contrast jumps ($\kappa_1 \neq \kappa_2$), the classical Delaunay sum condition $\alpha_1 + \alpha_2 \le \pi$ fails to ensure non-positivity.
Consequently, our standing geometric assumption requires $\mathcal T_h$ to be elementwise non-obtuse (all internal angles $\le \pi/2$ in $d=2$, or all dihedral angles $\le \pi/2$ in $d=3$), guaranteeing $\cot \ge 0$ individually and thus $A_{h,ij} \le 0$ for all $i \neq j$.

We emphasize, however, that global $M$-matrix inverse positivity does not guarantee that the interior partition sum $q_0(p) := \sum_{i\in\mathcal I_0} \chi_i^{\rm ms}(p)$ is strictly positive everywhere in $D$.
Indeed, whenever a boundary coarse element has all its vertices on $\partial D$ (such as coarse corner simplices), it lies entirely outside the support of all interior coarse basis functions $\{\chi_i^{\rm ms}\}_{i\in\mathcal I_0}$.
At any interior fine node $p$ within such a corner element, $q_0(p) = 0$ identically, precluding a global pointwise division $(\mathcal T_0 w)(p) = w(p)/q_0(p)$.
Furthermore, even when $q_0(p) > 0$, an arbitrary localized fine-grid trace $R_i(\mathcal T_0 w)$ does not generically belong to the chosen reduced snapshot subspace $W_i \subsetneq V_h(\omega_i)$ (e.g., in discrete harmonic snapshot spaces).

In standard GMsFEM implementations, this restriction is handled either by retaining boundary coarse patches equipped with snapshot spaces satisfying homogeneous Dirichlet boundary conditions on $\partial D \cap \omega_i$, or by establishing error estimates conditional on the target state possessing an interior-patch representation.
In strict alignment with \cref{thm:implementation-matched-infimal-multiscale}, we therefore formulate the analysis under the explicit hypothesis that the target fine-grid state $w = u_h(\theta)$ admits an admissible decomposition $\boldsymbol r = (r_i)_{i\in\mathcal I_0} \in \mathcal D_h(w)$ with $r_i \in W_i$ satisfying $w = \sum_{i\in\mathcal I_0} \mathcal M_i r_i$.
Under this hypothesis ($\mathcal D_h(w) \neq \emptyset$), the admissible set $\mathcal D_h(w)$ is a nonempty, closed affine subspace of the finite-dimensional product space $\prod_{i\in\mathcal I_0} W_i$.
Since the quadratic objective $\boldsymbol r \mapsto \sum_{i\in\mathcal I_0} \frac{C_{\times,i}^{\rm ref}}{\lambda_{i,L_i+1}} \mathfrak R_i(r_i)$ is continuous, non-negative, and convex, the infimum in \cref{eq:infimal-multiscale-functional} is attained by a minimizer, ensuring $\mathcal E_h(w) < \infty$.

\subsubsection{Step 3: Contrast-independent product bound and global assembly}
We next establish the uniform bound on $C_{\times,i}^{\rm ref}$. Let $T\in\mathcal T_h$ be a fine element (triangle in $d=2$ or tetrahedron in $d=3$) contained in $\omega_i$.
Because $\mathcal M_i z = \mathcal I_h(\chi_i^{\rm ms} z)$ is piecewise linear on $T$, its gradient $\nabla \mathcal M_i z$
is constant on $T$. Decomposing $z|_T$ into its element mean $\bar z_T$ and zero-mean fluctuation $\widetilde z_T := z|_T - \bar z_T$,
the nodal interpolation satisfies:
\[
  \nabla \mathcal I_h(\chi_i^{\rm ms} z)|_T
  =
  \bar z_T \nabla \chi_i^{\rm ms}|_T + \nabla \mathcal I_h(\chi_i^{\rm ms} \widetilde z_T)|_T.
\]
By standard finite element norm equivalence on the shape-regular reference element $\widehat T$ and scaling back to $T$,
there exists a geometric constant $C_1>0$ depending only on the shape regularity of $\mathcal T_h$ such that
\[
  \lVert\nabla \mathcal I_h(\chi_i^{\rm ms} \widetilde z_T)\rVert_{L^2(T)}^2
  \le
  C_1 \lVert\nabla z\rVert_{L^2(T)}^2,
  \qquad
  \lVert\bar z_T \nabla \chi_i^{\rm ms}\rVert_{L^2(T)}^2
  \le
  C_1 \lVert\nabla\chi_i^{\rm ms}\rVert_{L^\infty(T)}^2 \lVert z\rVert_{L^2(T)}^2.
\]
Applying the inequality $(a+b)^2 \le 2(a^2+b^2)$ yields:
\begin{equation}
  \lVert\nabla \mathcal I_h(\chi_i^{\rm ms} z)\rVert_{L^2(T)}^2
  \le
  2 C_1 \left( \lVert\nabla z\rVert_{L^2(T)}^2 + \lVert\nabla\chi_i^{\rm ms}\rVert_{L^\infty(T)}^2 \lVert z\rVert_{L^2(T)}^2 \right).
  \label{eq:app-triangle-bound}
\end{equation}
Multiplying \cref{eq:app-triangle-bound} by the constant value $\kappa_{\rm ref}|_T$, summing over all fine elements
$T\subset\omega_i$, and using the pointwise relation $|\nabla\chi_i^{\rm ms}|^2 \le \sum_{j\in\mathcal N_H} |\nabla\chi_j^{\rm ms}|^2$,
we obtain:
\begin{align*}
  \lVert\mathcal M_i z\rVert_{a_{\rm ref}}^2
  &=
  \int_{\omega_i} \kappa_{\rm ref} |\nabla \mathcal I_h(\chi_i^{\rm ms} z)|^2\,dx \\
  &\le
  2 C_1 \left( \int_{\omega_i} \kappa_{\rm ref} |\nabla z|^2\,dx + \int_{\omega_i} \kappa_{\rm ref} \sum_{j\in\mathcal N_H} |\nabla\chi_j^{\rm ms}|^2 |z|^2\,dx \right) \\
  &\le
  2 C_1 \left( a_i^{\rm ref}(z,z) + s_i^+(z,z) \right).
\end{align*}
Setting $C_{\rm sh} := 2 C_1$, we deduce that $C_{\times,i}^{\rm ref} \le C_{\rm sh}$ for all $i\in\mathcal I_0$.
Because $\kappa_{\rm ref}|_T$ is cellwise constant, it factors out identically in both terms of the local energy comparison. Thus, $C_{\rm sh}$ depends solely on the aspect ratios of the elements in $\mathcal T_h$ and the uniform gradient bound on $\{\chi_i^{\rm ms}\}$, remaining strictly independent of the mesh sizes $h, H$ and the physical contrast ratio $\kappa_{\max}/\kappa_{\min}$. This establishes the contrast-independent stability of the partition-of-unity interpolation constant in \cref{eq:uniform-nodal-product-bound}.

Finally, fix $w = u_h(\theta)$ and let $\boldsymbol r\in\mathcal D_h(u_h(\theta))$ be an arbitrary admissible decomposition.
Consider the multiscale approximant $v_{\rm ms} := \sum_{i\in\mathcal I_0} \mathcal M_i \Pi_i^{(L_i)} r_i \in V_{\rm ms,0}$.
Because at most $N_{\rm ov} \le 4$ coarse neighborhoods overlap at any point in $D$, Cauchy--Schwarz implies
$\lVert\sum_i \mathcal M_i e_i\rVert_{a_{\rm ref}}^2 \le N_{\rm ov} \sum_i \lVert\mathcal M_i e_i\rVert_{a_{\rm ref}}^2$.
Combining this with the local stability bound and spectral tail inequality \cref{eq:app-tail-sum} yields:
\[
  \lVert u_h(\theta) - v_{\rm ms}\rVert_{a_{\rm ref}}^2
  \le
  N_{\rm ov} \sum_{i\in\mathcal I_0} C_{\times,i}^{\rm ref} \left( a_i^{\rm ref}(e_i, e_i) + s_i^+(e_i, e_i) \right)
  \le
  N_{\rm ov} \sum_{i\in\mathcal I_0} \frac{C_{\times,i}^{\rm ref}}{\lambda_{i,L_i+1}} \mathfrak R_i(r_i).
\]
By the comparability definition \cref{eq:reference-comparability-ratio}, $\lVert v\rVert_{a_\theta}^2 \le \rho_+(\theta) \lVert v\rVert_{a_{\rm ref}}^2$.
Applying Galerkin best approximation for $u_G(\theta)$ in $V_{\rm ms,0}$ gives:
\[
  \lVert u_h(\theta) - u_G(\theta)\rVert_{a_\theta}^2
  \le
  \lVert u_h(\theta) - v_{\rm ms}\rVert_{a_\theta}^2
  \le
  \rho_+(\theta) N_{\rm ov} \sum_{i\in\mathcal I_0} \frac{C_{\times,i}^{\rm ref}}{\lambda_{i,L_i+1}} \mathfrak R_i(r_i).
\]
Taking the infimum over all $\boldsymbol r\in\mathcal D_h(u_h(\theta))$ completes the proof of \cref{eq:implementation-matched-infimal-estimate}.

We emphasize that while the partition-of-unity stability constant $C_{\rm sh}$ is strictly contrast-independent, \cref{eq:implementation-matched-infimal-estimate} bounds the squared energy error by an inverse eigenvalue factor $\mathcal O(\lambda_{\min,L+1}^{-1})$ (corresponding to an energy-norm decay of $\mathcal O(\lambda_{\min,L+1}^{-1/2})$ under a controlled numerator). The overall multiscale reduction error may retain contrast dependence through the parameter comparability ratio $\rho_+(\theta)$, the omitted eigenvalues $\lambda_{i,L_i+1}$, and the decomposition norms $\mathfrak R_i(r_i)$.
Full contrast robustness of the state error additionally requires that high-conductivity channels are isolated such that the omitted eigenvalues remain bounded away from zero, and that target states admit decompositions whose energies $\mathfrak R_i(r_i)$ are bounded uniformly with respect to the physical contrast ratio.
\qed

\subsection{Proof of Proposition~\ref{prop:finite-penalty-relaxation}}
\label{app:finite-penalty-relaxation-proof}

Recall the notation $B = H_f^\top W_y H_f + \lambda_u I_{N_c}\succ0$ and $h = H_f^\top W_y y$
defined in \cref{eq:observation-regularization-block}. Abbreviating $A := A_f(\theta)$,
the normal equations \cref{eq:projected-state-normal-equations} for the projected state
$c_\varepsilon(\theta)$ with weak weight matrix $W_W(\varepsilon) = \varepsilon^{-2}\overline W_W$ read:
\begin{equation*}
  \left(\frac{A^\top \overline W_W A}{\varepsilon^2} + B\right) c_\varepsilon(\theta)
  =
  \frac{A^\top \overline W_W b_f}{\varepsilon^2} + h.
\end{equation*}
Multiplying both sides by $\varepsilon^2$ yields the equivalent nonsingular system:
\begin{equation}
  \left(A^\top \overline W_W A + \varepsilon^2 B\right) c_\varepsilon(\theta)
  =
  A^\top \overline W_W b_f + \varepsilon^2 h.
  \label{eq:appendix-scaled-finite-penalty-identity}
\end{equation}
Because $c_G(\theta)$ satisfies the exact Galerkin balance $A c_G(\theta) = b_f$, we have
$A^\top \overline W_W b_f = A^\top \overline W_W A c_G(\theta)$. Adding $\varepsilon^2 B c_G(\theta)$
to both sides gives:
\begin{equation}
  \left(A^\top \overline W_W A + \varepsilon^2 B\right) c_G(\theta)
  =
  A^\top \overline W_W b_f + \varepsilon^2 B c_G(\theta).
  \label{eq:appendix-hard-state-in-scaled-system}
\end{equation}
Subtracting \cref{eq:appendix-hard-state-in-scaled-system} from
\cref{eq:appendix-scaled-finite-penalty-identity} leads to:
\begin{equation}
  \left(A^\top \overline W_W A + \varepsilon^2 B\right) \left(c_\varepsilon(\theta) - c_G(\theta)\right)
  =
  \varepsilon^2 \left(h - B c_G(\theta)\right).
  \label{eq:appendix-error-system}
\end{equation}
By the uniform ellipticity assumption \ref{ass:coefficients}, $A = A_f(\theta)$ is nonsingular.
Since $\overline W_W\succ0$ and $B\succ0$, the matrix $A^\top \overline W_W A + \varepsilon^2 B$ is strictly
positive definite for all $\varepsilon > 0$. Inverting this matrix in \cref{eq:appendix-error-system}
proves the exact displacement identity \cref{eq:exact-finite-penalty-displacement}.

To establish the norm estimate, observe that for any $\varepsilon > 0$,
\[
  A^\top \overline W_W A + \varepsilon^2 B \succeq A^\top \overline W_W A.
\]
Consequently, the spectral norm of its inverse is bounded by:
\[
  \left\lVert\left(A^\top \overline W_W A + \varepsilon^2 B\right)^{-1}\right\rVert_2
  \le
  \left\lVert\left(A^\top \overline W_W A\right)^{-1}\right\rVert_2
  =
  \frac{1}{\sigma_{\min}\big(\overline W_W^{1/2} A\big)^2}.
\]
Taking Euclidean norms in \cref{eq:exact-finite-penalty-displacement} yields:
\[
  \lVert c_\varepsilon(\theta) - c_G(\theta)\rVert_2
  \le
  \frac{\varepsilon^2}{\sigma_{\min}\big(\overline W_W^{1/2} A_f(\theta)\big)^2}
  \lVert h - B c_G(\theta)\rVert_2.
\]
Applying the linear reconstruction mapping $\lVert\Phi_f c\rVert_X \le \lVert\Phi_f\rVert_{\ell^2\to X} \lVert c\rVert_2$
completes the proof of \cref{eq:conditioning-explicit-finite-penalty-bound}.
\qed


\section{Projected-State Sensitivity and Lipschitz Continuity}
\label{app:conditional-inverse-optimization}

This appendix provides the explicit Jacobian formulation of the projected-state mapping
$\theta\mapsto c_\star(\theta)$ and proves \cref{prop:projected-state-lipschitz-reductions}.

\subsection{Proof of Proposition~\ref{prop:projected-state-lipschitz-reductions}}
\label{app:projected-state-lipschitz-proof}

Recall the symmetric positive definite normal Gram matrix and load vector from
\cref{eq:normal-equations-rep}:
\begin{equation}
  M(\theta) := H_f^\top W_y H_f + A_f(\theta)^\top W_W A_f(\theta) + \lambda_u I_{N_c},
  \qquad
  g(\theta) := H_f^\top W_y y + A_f(\theta)^\top W_W b_f.
  \label{eq:app-normal-system}
\end{equation}
The unique projected state satisfies $M(\theta) c_\star(\theta) = g(\theta)$ with
$M(\theta) \succeq \lambda_u I_{N_c}$. By assumption \ref{ass:parameter-regularity},
$\theta\mapsto A_f(\theta)$ is continuously differentiable on an open neighborhood of $\Theta$.
Differentiating the identity $M(\theta) c_\star(\theta) = g(\theta)$ with respect to coordinate
$\theta_r$ ($r=1,\ldots,N_\theta$) yields:
\[
  \frac{\partial M(\theta)}{\partial\theta_r} c_\star(\theta)
  +
  M(\theta) \frac{\partial c_\star(\theta)}{\partial\theta_r}
  =
  \frac{\partial g(\theta)}{\partial\theta_r}.
\]
Solving for the partial derivative of $c_\star(\theta)$ gives:
\begin{equation}
  \frac{\partial c_\star(\theta)}{\partial\theta_r}
  =
  M(\theta)^{-1}
  \left(
    \frac{\partial g(\theta)}{\partial\theta_r}
    -
    \frac{\partial M(\theta)}{\partial\theta_r} c_\star(\theta)
  \right).
  \label{eq:app-state-derivative-raw}
\end{equation}
Evaluating the parameter derivatives of $g(\theta)$ and $M(\theta)$ directly from \cref{eq:app-normal-system}:
\begin{align}
  \frac{\partial g(\theta)}{\partial\theta_r}
  &=
  \left(\frac{\partial A_f(\theta)}{\partial\theta_r}\right)^\top W_W b_f,
  \label{eq:app-grad-g-clean} \\
  \frac{\partial M(\theta)}{\partial\theta_r}
  &=
  \left(\frac{\partial A_f(\theta)}{\partial\theta_r}\right)^\top W_W A_f(\theta)
  +
  A_f(\theta)^\top W_W \frac{\partial A_f(\theta)}{\partial\theta_r}.
  \notag
\end{align}
Substituting these expressions into \cref{eq:app-state-derivative-raw} and grouping terms
via the discrete weak residual $r_W(c_\star(\theta),\theta) = A_f(\theta)c_\star(\theta) - b_f$,
we obtain the compact derivative formula:
\begin{equation}
  \frac{\partial c_\star(\theta)}{\partial\theta_r}
  =
  -M(\theta)^{-1}
  \left[
    \left(\frac{\partial A_f(\theta)}{\partial\theta_r}\right)^\top W_W r_W(c_\star(\theta),\theta)
    +
    A_f(\theta)^\top W_W \frac{\partial A_f(\theta)}{\partial\theta_r} c_\star(\theta)
  \right].
  \label{eq:app-state-derivative-compact}
\end{equation}

Because $\theta\mapsto A_f(\theta)$ is $C^1$ and $M(\theta)\succeq\lambda_u I_{N_c}$ with $\lambda_u>0$,
the Jacobian matrix $D_\theta c_\star(\theta) \in \mathbb R^{N_c\times N_\theta}$ is continuous on $\Theta$.
Under the explicit hypothesis of Proposition~\ref{prop:projected-state-lipschitz-reductions}, $\mathcal N \subset \Theta$ is a prescribed compact convex subset containing $\theta^\dagger$, $\theta_\star$, and $\widehat\theta$ (which exists, for example, by taking $\mathcal N$ as the closed convex hull of a compact sublevel set of $\varphi$ when $\Theta$ is convex).
By the extreme value theorem, the spectral norm of $D_\theta c_\star(\theta)$ achieves a finite maximum on $\mathcal N$,
guaranteeing that
\[
  L_c(\mathcal N) = \sup_{\theta\in\mathcal N} \lVert D_\theta c_\star(\theta)\rVert_{2\to2} < \infty.
\]

Because $\mathcal N \subset \Theta$ is convex, for any two parameter points $\theta, \mu\in\mathcal N$, the line segment $\gamma(t) := \mu + t(\theta - \mu)$ for $t\in[0,1]$ remains entirely contained within $\mathcal N \subset \Theta$. This ensures that the integrand $D_\theta c_\star(\gamma(t))$ is well-defined and uniformly bounded along the entire integration path.
Applying the fundamental theorem of calculus:
\[
  c_\star(\theta) - c_\star(\mu)
  =
  \int_0^1 D_\theta c_\star\big(\mu + t(\theta - \mu)\big) (\theta - \mu)\,dt.
\]
Taking the Euclidean norm and applying the integral triangle inequality yields:
\[
  \lVert c_\star(\theta) - c_\star(\mu)\rVert_2
  \le
  \int_0^1 \big\lVert D_\theta c_\star\big(\mu + t(\theta - \mu)\big)\big\rVert_{2\to2} \lVert\theta - \mu\rVert_2\,dt
  \le
  L_c(\mathcal N) \lVert\theta - \mu\rVert_2.
\]
Multiplying by the operator norm of the linear state synthesis operator $\Phi_f$ confirms
\cref{eq:compact-projected-state-lipschitz}. Finally, specializing $(\theta,\mu)$ to the pairs
$(\theta^\dagger, \widehat\theta)$, $(\theta^\dagger, \theta_\star)$, and $(\theta_\star, \widehat\theta)$
directly yields the reductions in \cref{eq:inverse-optimization-state-reductions}.
\qed

\bibliographystyle{unsrtnat-nodoi}
\bibliography{references}

\end{document}